\documentclass[11pt]{article}

\usepackage{amscd,amstex,amssymb}

\def\eqref#1{(\ref{#1})}
\newcommand{\goth}{\mathfrak}

\newcommand{\arrow}{{\:\longrightarrow\:}}
\newcommand{\Z}{{\Bbb Z}}
\newcommand{\C}{{\Bbb C}}
\newcommand{\R}{{\Bbb R}}

\newcommand{\1}{\sqrt{-1}\:}
\newcommand{\restrict}[1]{{\left|_{{\phantom{|}\!\!}_{#1}}\right.}}

\renewcommand{\c}[1]{{\cal #1}}
\newcommand{\calo}{{\cal O}}

\edef\newboxtimes{\mathop{\boxtimes}}
\renewcommand{\tilde}{\widetilde}
\renewcommand{\bar}{\overline}
\renewcommand{\phi}{\varphi}
\renewcommand{\epsilon}{\varepsilon}
\renewcommand{\geq}{\geqslant}
\renewcommand{\leq}{\leqslant}

\newcommand{\im}{\operatorname{im}}
\newcommand{\End}{\operatorname{End}}
\newcommand{\Id}{\operatorname{Id}}

\newcommand{\Aut}{\operatorname{Aut}}
\newcommand{\codim}{\operatorname{codim}}
\newcommand{\slope}{\operatorname{slope}}

\newcommand{\Def}{\operatorname{Def}}

\newcommand{\comment}[1]{{}}

\def\blacksquare{\hbox{\vrule width 4pt height 4pt depth 0pt}}
\def\endproof{\blacksquare}

\makeatletter

\@ifundefined{Bbb}
     {\newcommand{\Bbb}[1]{{\mathbb #1}}}%
{}

\newcommand{\ps@verbit}{%
  \renewcommand{\@oddhead}{%
          \scriptsize
          {Subvarieties of the Hilbert scheme}
          \hfil\tiny {M. Verbitsky, \ \ \ April 5, 1997, 
			revised April 25, Oct 31}}
  \renewcommand{\@evenhead}{\@oddhead}
  \renewcommand{\@oddfoot}{\hfil\thepage\hfil}
  \renewcommand{\@evenfoot}{\@oddfoot}}
 
\newcounter{Mycounter}[section]
\newcounter{lemma}[section]
\renewcommand{\thelemma}{{Lemma \thesection.\arabic{lemma}}}
\newcommand{\lemma}{%
     \setcounter{lemma}{\value{Mycounter}}
     \refstepcounter{lemma}
     \stepcounter{Mycounter}
     {\bf \thelemma:\ }}

\newcounter{claim}[section]
\renewcommand{\theclaim}{{Claim \thesection.\arabic{claim}}}
\newcommand{\claim}{%
     \setcounter{claim}{\value{Mycounter}}
     \refstepcounter{claim}
     \stepcounter{Mycounter}
     {\bf \theclaim:\ }}

\newcounter{sublemma}[section]
\renewcommand{\thesublemma}{{Sublemma \thesection.\arabic{sublemma}}}
\newcommand{\sublemma}{%
     \setcounter{sublemma}{\value{Mycounter}}
     \refstepcounter{sublemma}
     \stepcounter{Mycounter}
     {\bf \thesublemma:\ }}

\newcounter{corollary}[section]
\renewcommand{\thecorollary}{{Corollary \thesection.\arabic{corollary}}}
\newcommand{\corollary}{%
     \setcounter{corollary}{\value{Mycounter}}
     \refstepcounter{corollary}
     \stepcounter{Mycounter}
     {\bf \thecorollary:\ }}

\newcounter{theorem}[section]
\renewcommand{\thetheorem}{{Theorem \thesection.\arabic{theorem}}}
\newcommand{\theorem}{%
     \setcounter{theorem}{\value{Mycounter}}
     \refstepcounter{theorem}
     \stepcounter{Mycounter}
     {\bf \thetheorem:\ }}

\newcounter{conjecture}[section]
\newcounter{proposition}[section]
\renewcommand{\theproposition}
       {{Proposition \thesection.\arabic{proposition}}}
\newcommand{\proposition}{%
     \setcounter{proposition}{\value{Mycounter}}
     \refstepcounter{proposition}
     \stepcounter{Mycounter}
     {\bf \theproposition:\ }}

\newcounter{definition}[section]
\renewcommand{\thedefinition}
       {{Definition \thesection.\arabic{definition}}}
\newcommand{\definition}{%
     \setcounter{definition}{\value{Mycounter}}
     \refstepcounter{definition}
     \stepcounter{Mycounter}
     {\bf \thedefinition:\ }}

\newcounter{example}[section]
\renewcommand{\theexample}{{Example \thesection.\arabic{example}}}
\newcommand{\example}{%
     \setcounter{example}{\value{Mycounter}}
     \refstepcounter{example}
     \stepcounter{Mycounter}
     {\bf \theexample:\ }}

\newcounter{remark}[section]
\renewcommand{\theremark}{{Remark \thesection.\arabic{remark}}}
\newcommand{\remark}{%
     \setcounter{remark}{\value{Mycounter}}
     \refstepcounter{remark}
     \stepcounter{Mycounter}
     {\bf \theremark:\ }}

\newcounter{problem}[section]
\newcounter{question}[section]
\@addtoreset{equation}{section}
\@addtoreset{footnote}{section}
\makeatother

\begin{document}

\begin{center}
{\Large\bf
Trianalytic subvarieties of the  
Hilbert scheme of points on a K3 surface}\\[4mm]
Misha Verbitsky,\footnote{Supported by the NSF grant 9304580}\\[4mm]
{\tt verbit@@thelema.dnttm.rssi.ru, verbit@@ihes.fr}
\end{center}

\hfill

{\small 
\hspace{0.2\linewidth}
\begin{minipage}[t]{0.7\linewidth}
Let $X$ be a hyperk\"ahler manifold . Trianalytic subvarieties of 
$X$ are subvarieties which are complex analytic
with respect to all complex structures
induced by the hyperk\"ahler structure.
Given a K3 surface $M$, the Hilbert scheme classifying 
zero-dimensional subschemes of $M$ admits a hyperk\"ahler
structure. We show that for
$M$ generic, there are no 
trianalytic subvarieties of the Hilbert scheme.
This implies that a generic deformation of the Hilbert
scheme of K3 has no proper complex subvarieties. 
\end{minipage}
}

\tableofcontents


\section{Introduction}


\hfill

For basic results and definitions of
hyperk\"ahler geometry, see
\cite{_Besse:Einst_Manifo_}. 

\hfill

This Introduction is independent from the rest of this paper.

\subsection{An overview}

An {\bf almost hypercomplex} manifold $M$ is a 
manifold equipped with an action of the 
quaternion algebra ${\Bbb H}$ on its tangent bundle. 
The manifold $M$ is called {\bf hypercomplex} if for all algebra
embedding $\C \stackrel{\iota}\hookrightarrow \Bbb H$,
the corresponding almost complex structure $I_\iota$
is integrable. A manifold $M$ is called {\bf hyperk\"ahler}
if, on top of that, $M$ is equipped with a Riemannian metric
which is K\"ahler with respect to the
complex structures $I_\iota$, for all 
embeddings $\C \stackrel{\iota}\hookrightarrow \Bbb H$.
The complex structures $I_\iota$ are called 
{\bf induced complex structures}; the corresponding
K\"ahler manifold is denoted by $(M, I_\iota)$.

For a more formal definition of a hyperk\"ahler manifold, see 
\ref{_hyperkahler_manifold_Definition_}.
The notion of a hyperk\"ahler manifold was introduced by E. Calabi,
\cite{_Calabi_}.

Clearly, the real dimension of $M$ is divisible by 4.
For $\dim_\R M= 4$, there are only two classes of compact 
hyperk\"ahler manifolds: compact tori and K3 surfaces. 

Let $M$ be a complex surface and $M^{(n)}$ be its $n$-th 
symmetric power, $M^{(n)} = M^n/S_n$. The variety $M^{(n)}$
admits a natural desingularization  $M^{[n]}$, called
{\bf Hilbert scheme of points}, or {\bf Hilbert scheme} 
for short. For its construction and additional results,
see the lectures of H. Nakajima, \cite{_Nakajima_}.

Most importantly, $M^{[n]}$ admits a hyperk\"ahler metrics
whenever the surface $M$ is compact and hyperk\"ahler (\cite{_Beauville_}).
This way, Beauville constructed 
two series of examples of hyperk\"ahler manifolds,
associated with a torus
\footnote{There is a natural torus action on its Hilbert scheme;
to obtain the Beauville's hyperk\"ahler manifold, we must take
the quotient by this action.}
 and a K3 surface. It was conjectured that
all  hyperk\"ahler manifolds $X$ with $H^1(X) =0$, $H^{2,0}(X)=\C$
are deformationally equivalent to one of these examples. 
We study the complex and hyperk\"ahler
geometry of $M^{[n]}$ for $M$ a ``sufficiently generic'' 
K3 surface, in order to construct counterexamples to 
this conjecture.

\hfill

Let $M$ be a hyperk\"ahler manifold. A 
{\bf trianalytic subvariety} of $M$ is a closed subset
which is complex analytic with respect to any of the induced complex
structures. It was proven in \cite{_Verbitsky:Symplectic_II_} 
that for all induced complex structures $I$, except maybe a
countable number, all complex subvarieties of $(M,I)$ are
trianalytic (see also \ref{_generic_are_dense_Proposition_}). 
This reduces the study of complex subvarieties 
of ``sufficiently generic'' 
deformations\footnote{By {\bf deformations of $M$} we understand complex
manifolds which are deformationally equivalent to $M$.}
of $M$ to the study of trianalytic subvarieties.

Trianalytic subvarieties of hyperk\"ahler manifolds were
studied at length in
\cite{_Verbitsky:Symplectic_II_} and \cite{_Verbitsky:DesinguII_}.
Of the results obtained in this study, the most important ones are  
 Desingularization Theorem (\ref{_desingu_Theorem_}) 
and the cohomological criterion of trianalyticity
(\ref{_G_M_invariant_implies_trianalytic_Theorem_}). 
The aim of the present paper is to 
obtain the following theorem.

\hfill

\theorem \label{_no_triana_subva_Introdu_Theorem_}
Let $M$ be a complex K3 surface without automorphisms, 
$\c H$ a hyperk\"ahler structure on $M$,
and $I$ an induced complex structure on $M$ which is Mumford-Tate 
generic in the class of induced complex structures.%
\footnote{For the definition of Mumford-Tate generic,
see \ref{_generic_manifolds_Definition_}.}
Let $M^{[n]}$ be a Hilbert scheme of points on  
the complex surface $(M, I)$. Pick any hyperk\"ahler
structure on $M^{[n]}$, compatible with the complex structure.
Then, $M^{[n]}$ has no proper trianalytic subvarieties.

{\bf Proof:} This is \ref{_no_triana_subva_of_Hilb_Theorem_}.
\endproof

\hfill

In the forthcoming paper, we construct a 
21-dimensional family of
compact hyperk\"ahler manifolds $M_x$, with
\begin{equation}\label{_non-decompo_Equation_}
  H^{1}(M) =0,  \ \  H^{2,0}(M)=\C
\end{equation}
which {\it have} proper trianalytic subvarieties. This 
leads to assertion that these manifolds are not 
deformations of $M^{[n]}$. 

As another application of 
\ref{_no_triana_subva_Introdu_Theorem_},
we obtain that a generic complex deformation of
$M^{[n]}$ has no proper closed complex 
subvarieties (\ref{_no_comple_subva_of_gen_Hilb_Corollary_}).

A version of this result is true for a compact torus.
For a generic complex structure $I$ on a complex torus
$T$, the complex manifold $(T,I)$ has no proper subvarieties.
This is easy to see from the fact that
the group $H^{p,p}(T,I) \cap H^{2p}(T, \Z)$ is empty.
For a Hilbert scheme of K3, this Hodge-theoretic
argument does not work. In fact, there are integer
cycles in $H^{2p,2p}(M^{[n]}, I)$ for all complex structures
on $M^{[n]}$. As \ref{_no_comple_subva_of_gen_Hilb_Corollary_}
implies, these cycles cannot be represented by subvarieties.
This gives a counterexample to the Hodge conjecture.

Of course, for a generic complex structure $I$, the
manifold $(M^{[n]}, I)$ is not algebraic. There are many other
counterexamples to the Hodge conjecture for non-algebraic
manifolds.

\hfill

In our approach to the study of trianalytic subvarieties of the
Hilbert scheme, we introduce the concept of the {\bf universal
subvariety} of the Hilbert scheme (\ref{_Unive_subva_Definition_}).
For a complex surface $M$, the local automorphisms
$\gamma:\; U \arrow U$ of $M\supset U$ act on the corresponding
open subsets $U^{[n]}\subset M^{[n]}$ of the Hilbert scheme.
A universal subvariety of $M^{[n]}$ is a subvariety which is preserved by
all automorphisms obtained this way (see \ref{_Unive_subva_Definition_}
for a more precise statement). 

We show that a trianalytic subvariety of a Hilbert scheme of
a K3 surface $M$ is universal, assuming that $M$ is Mumford-Tate
generic with respect to some hyperk\"ahler structure
and has no complex automorphisms 
(\ref{_triana_subva_universal_Theorem_}).

\subsection{Contents}

The paper is organized as follows.

\begin{itemize} 

\item Section \ref{_basics_hyperka_Section_}
is taken with preliminary conventions and basic theorems. 
We define hyperk\"ahler manidols and formulate Yau's theorem
on the existence of hyperk\"ahler structures on compact
holomorphically symplectic manifolds of K\"ahler type.
Furthermore, we define trianalytic subvarieties of 
hyperk\"ahler manifolds and recall the basic properties of trianalytic 
subvarieties, following \cite{_Verbitsky:Symplectic_II_} 
and \cite{_Verbitsky:DesinguII_}. There are no new
results in Section \ref{_basics_hyperka_Section_},
and nothing unknown to the reader acquainted with 
the literature. 

\item In Section \ref{_appendix_subva_produ_Section_},
we apply the desingularization theorem to the subvarieties of 
$M^n$, where $M$ is a generic K3 surface.
We classify trianalytic subvarieties of $M^n$
and describe them explicitly. This section is
independent from the rest of this paper.

\item We study the Hilbert scheme $M^{[n]}$ of a smooth holomorphic
symplectic complex surface $M$
in Section \ref{_Hilbert_sche_Section_}. We give its definition
and explain the construction of the holomorphic symplectic
structure on $M^{[n]}$. By Yau's proof of Calabi conjecture
(\ref{_symplectic_=>_hyperkahler_Proposition_}),
this implies that $M^{[n]}$ admits a hyperk\"ahler structure.

Using perverse sheaves, we write down
the cohomology of $M^{[n]}$ in terms of diagonals in the symmetric
power $M^{(n)}$. This is done using the fact that the standard projection
$\pi:\; M^{[n]}\arrow M^{(n)}$ is a semi-small resolution of
the symmetric power $M^{(n)}$.
These results are well known (\cite{_Nakajima_}).

Further on, we apply the same type arguments to the
trianalytic subvarieties $X\subset M^{[n]}$. 
Using the holomorphic symplectic geometry, we show that
the map $\tilde X\stackrel {\pi\circ n}\arrow \pi(X)$ 
is a semi-small resolution, where 
$\tilde X\stackrel {n}\arrow X$ is the hyperk\"ahler
desingularization of $X$ (\ref{_desingu_Theorem_}). 
This gives an expression for the cohomology of
$\tilde X$. We don't use this result anywhere 
in this paper.

\item Section \ref{_Hilbert_sche_Section_} is 
heavily based on perverse sheaves (\cite{_Asterisque_100_}), 
and does not use  results of hyperk\"ahler geometry, except 
Desingularization Theorem (\ref{_desingu_Theorem_}). 

\item The following four sections
(Sections \ref{_unive_subva_Section_}--\ref{_triana_unive_subva_Section_})
are dedicated to the study of universal properties of the
Hilbert scheme. 

\begin{itemize}

\item In Section \ref{_unive_subva_Section_}, we give a definition
of a universal subvariety of a Hilbert scheme. A
relative dimension of a universal subvariety is a 
dimension of the generic fiber of the
projection $\pi:\; X \arrow \pi(X)$, where
$\pi:\; M^{[n]}\arrow M^{(n)}$ is the standard
projection of the Hilbert scheme to the symmetric
power of $M$. We classify and describe explicitly 
the universal subvarieties of relative dimension $0$.
Results of Section \ref{_unive_subva_Section_} are in no
way related to the hyperk\"ahler geometry. 

\item In Section \ref{_subva_gene_fini_Section_},
we study the Hilbert scheme of a K3 surface $M$,
assuming that $M$ is Mumford-Tate generic with
respect to some hyperk\"ahler structure.
We consider subvarieties $X\subset M^{[n]}$,
such that $X$ is projected generically finite
to $\pi(X)\subset M$, and $\pi(X)$ is 
a diagonal in $M^{(n)}$.
We use the theory of Yang-Mills connections
and Uhlenbeck--Yau theorem, in order to show that
such subvarieties are universal, in the sense of
Section \ref{_unive_subva_Section_}.

\item Section \ref{_specia_subva_Section_}
is completely parallel to Section \ref{_unive_subva_Section_}.
We define {\bf special subvarieties} of the Hilbert scheme,
which are similar to the universal subvarieties,
with some conditions relaxed. Whereas 
universal subvarieties are subvarieties which are
fixed by all local automorphisms of $M^{[n]}$ coming from
$M$, special subvarieties are the subvarieties 
fixed by all the local automorphisms coming from
$M$ which preserve a finite subset $S\subset M$.
As in Section \ref{_unive_subva_Section_},
we classify and describe explicitely the
special subvarieties of relative dimension zero. 
Using results of Section \ref{_subva_gene_fini_Section_},
we study the subvarieties $X\subset M^{[n]}$,
such that $X$ is projected generically finite
to $\pi(X)$, where $M$ is a generic K3 surface. 
We show that all such subvarieties are 
special of relative dimension zero. 

\item In Section \ref{_triana_unive_subva_Section_},
we study the deformations of subvarieties of 
the Hilbert scheme of a K3 surface. 
The deformations of special subvarieties of relative dimension zero
are easy to study using the explicit description given in
Section \ref{_specia_subva_Section_}. The
deformations of trianalytic subvarieties were studied at length in
\cite{_Verbitsky:Deforma_}. In conjunction, these
results lead to the assertion that all trianalytic
subvarieties of $M^{[n]}$ are universal, in the sense
of Section \ref{_unive_subva_Section_}.

\end{itemize}

\item In Section \ref{_last_Section_}, we study the second cohomology
of universal subvarieties $\c X_\alpha\stackrel \phi\hookrightarrow M^{[n]}$ 
of the Hilbert scheme $M^{[n]}$,
in assumption that $\c X_\alpha$ is trianalytic.
First of all, we show that $\c X_\alpha$ is birationally
equivalent to a hyperka\"hler manifold which is a product
of Hilbert schemes of $M$. Using Mukai's theorem, which states
that second cohomology of hyperka\"hler manifolds is a birational
invariant, we obtain a structure theorem for $H^2(\c X_\alpha)$.
Assuming that $\c X_\alpha$ is not a product of
hyperka\"hler manifolds, we show that the pullback
map $\phi^*:\; H^2(M^{[n]}) \arrow H^2(\c X_\alpha)$
is an isomorphism, and compute this map explicitly.
The second cohomology of a hyperka\"hler manifold 
$X$ is equipped with a canonical non-degenerate 
quadratic form $(\cdot,\cdot)_{\c B}$, defined up to 
a constant multiplier. This form is invariant under the
natural $SU(2)$-action on $H^2(X)$. We compute the pullback
of the form $(\cdot,\cdot)_{\c B}$ under the map  
$\phi^*:\; H^2(M^{[n]}) \arrow H^2(\c X_\alpha)$,
and show that it cannot be $SU(2)$-invariant.
Thus, $\phi^*$ is not compatible with the $SU(2)$-action 
on the second cohomology. This implies that $\phi$ cannot
be compatible with the hyperka\"hler structures on
$\c X_\alpha$, $M^{[n]}$. Therefore,
$M^{[n]}$ contains no trianalytic subvarieties.

\end{itemize}


\section{Hyperk\"ahler manifolds}
\label{_basics_hyperka_Section_}


\subsection{Hyperk\"ahler manifolds}

This subsection contains a compression of 
the basic and best known results 
and definitions from hyperk\"ahler geometry, found, for instance, in
\cite{_Besse:Einst_Manifo_} or in \cite{_Beauville_}.

\hfill

\definition \label{_hyperkahler_manifold_Definition_} 
(\cite{_Besse:Einst_Manifo_}) A {\bf hyperk\"ahler manifold} is a
Riemannian manifold $M$ endowed with three complex structures $I$, $J$
and $K$, such that the following holds.
 
\begin{description}
\item[(i)]  the metric on $M$ is K\"ahler with respect to these complex 
structures and
 
\item[(ii)] $I$, $J$ and $K$, considered as  endomorphisms
of a real tangent bundle, satisfy the relation 
$I\circ J=-J\circ I = K$.
\end{description}

\hfill 

The notion of a hyperk\"ahler manifold was 
introduced by E. Calabi (\cite{_Calabi_}).

\hfill

Clearly, a hyperk\"ahler manifold has a natural action of
the quaternion algebra ${\Bbb H}$ in its real tangent bundle $TM$. 
Therefore its complex dimension is even.
For each quaternion $L\in \Bbb H$, $L^2=-1$,
the corresponding automorphism of $TM$ is an almost complex
structure. It is easy to check that this almost 
complex structure is integrable (\cite{_Besse:Einst_Manifo_}).

\hfill

\definition \label{_indu_comple_str_Definition_} 
Let $M$ be a hyperk\"ahler manifold, $L$ a quaternion satisfying
$L^2=-1$. The corresponding complex structure on $M$ is called
{\bf an induced complex structure}. The $M$ considered as a K\"ahler
manifold is denoted by $(M, L)$. In this case,
the hyperk\"ahler structure is called {\bf combatible
with the complex structure $L$}.

Let $M$ be a compact complex variety. We say
that $M$ is {\bf of hyperk\"ahler type}
if $M$ admits a hyperk\"ahler structure
compatible with the complex structure.

\hfill

\hfill

\definition \label{_holomorphi_symple_Definition_} 
Let $M$ be a complex manifold and $\Theta$ a closed 
holomorphic 2-form over $M$ such that 
$\Theta^n=\Theta\wedge\Theta\wedge...$, is
a nowhere degenerate section of a canonical class of $M$
($2n=dim_\C(M)$).
Then $M$ is called {\bf holomorphically symplectic}.

\hfill

Let $M$ be a hyperk\"ahler manifold; denote the
Riemannian form on $M$ by $<\cdot,\cdot>$.
Let the form $\omega_I := <I(\cdot),\cdot>$ be the usual K\"ahler
form  which is closed and parallel
(with respect to the Levi-Civitta connection). Analogously defined 
forms $\omega_J$ and $\omega_K$ are
also closed and parallel. 
 
A simple linear algebraic
consideration (\cite{_Besse:Einst_Manifo_}) shows that the form
$\Theta:=\omega_J+\sqrt{-1}\omega_K$ is of
type $(2,0)$ and, being closed, this form is also holomorphic.
Also, the form $\Theta$ is nowhere degenerate, as another linear 
algebraic argument shows.
It is called {\bf the canonical holomorphic symplectic form
of a manifold M}. Thus, for each hyperk\"ahler manifold $M$,
and an induced complex structure $L$, the underlying complex manifold
$(M,L)$ is holomorphically symplectic. The converse assertion
is also true:

\hfill

\theorem \label{_symplectic_=>_hyperkahler_Proposition_}
(\cite{_Beauville_}, \cite{_Besse:Einst_Manifo_})
Let $M$ be a compact holomorphically
symplectic K\"ahler manifold with the holomorphic symplectic form
$\Theta$, a K\"ahler class 
$[\omega]\in H^{1,1}(M)$ and a complex structure $I$. 
Let $n=\dim_\C M$. Assume that
 $\int_M \omega^n = \int_M (Re \Theta)^n$.
Then there is a unique hyperk\"ahler structure $(I,J,K,(\cdot,\cdot))$
over $M$ such that the cohomology class of the symplectic form
$\omega_I=(\cdot,I\cdot)$ is equal to $[\omega]$ and the
canonical symplectic form $\omega_J+\1\omega_K$ is
equal to $\Theta$.

\hfill

\ref{_symplectic_=>_hyperkahler_Proposition_} immediately
follows from the conjecture of Calabi, pro\-ven by
Yau (\cite{_Yau:Calabi-Yau_}). 
\endproof

\hfill

Let $M$ be a hyperk\"ahler manifold. We identify the group $SU(2)$
with the group of unitary quaternions. This gives a canonical 
action of $SU(2)$ on the tangent bundle, and all its tensor
powers. In particular, we obtain a natural action of $SU(2)$
on the bundle of differential forms. 

\hfill

\lemma \label{_SU(2)_commu_Laplace_Lemma_}
The action of $SU(2)$ on differential forms commutes
with the Laplacian.
 
{\bf Proof:} This is Proposition 1.1
of \cite{_Verbitsky:Symplectic_I_}. \endproof
 
Thus, for compact $M$, we may speak of the natural action of
$SU(2)$ in cohomology.

\hfill

Further in this article, we use the following statement.

\hfill

\lemma \label{_SU(2)_inva_type_p,p_Lemma_} 
Let $\omega$ be a differential form over
a hyperk\"ahler manifold $M$. The form $\omega$ is $SU(2)$-invariant
if and only if it is of Hodge type $(p,p)$ with respect to all 
induced complex structures on $M$.

{\bf Proof:} This is \cite{_Verbitsky:Hyperholo_bundles_}, 
Proposition 1.2. \endproof

\subsection{Trianalytic subvarieties in compact hyperk\"ahler
manifolds.}
 
In this subsection, we give a definition and a few basic properties
of trianalytic subvarieties of hyperk\"ahler manifolds. 
We follow \cite{_Verbitsky:Symplectic_II_}.
 
\hfill
 
Let $M$ be a compact hyperk\"ahler manifold, $\dim_\R M =2m$.
 
\hfill
 
\definition\label{_trianalytic_Definition_} 
Let $N\subset M$ be a closed subset of $M$. Then $N$ is
called {\bf trianalytic} if $N$ is a complex analytic subset 
of $(M,L)$ for any induced complex structure $L$.
 
\hfill

\hfill
 
Let $I$ be an induced complex structure on $M$,
and $N\subset(M,I)$ be
a closed analytic subvariety of $(M,I)$, $dim_\C N= n$.
Consider the homology class 
represented by $N$. Let $[N]\in H^{2m-2n}(M)$ denote 
the Poincare dual cohomology class. Recall that
the hyperk\"ahler structure induces the action of 
the group $SU(2)$ on the space $H^{2m-2n}(M)$.
 
\hfill
 
\theorem\label{_G_M_invariant_implies_trianalytic_Theorem_} 
Assume that $[N]\in  H^{2m-2n}(M)$ is invariant with respect
to the action of $SU(2)$ on $H^{2m-2n}(M)$. Then $N$ is trianalytic.
 
{\bf Proof:} This is Theorem 4.1 of 
\cite{_Verbitsky:Symplectic_II_}.
\endproof
 
\remark \label{_triana_dim_div_4_Remark_}
Trianalytic subvarieties have an action of quaternion algebra in
the tangent bundle. In particular,
the real dimension of such subvarieties is divisible by 4.

\hfill

\definition \label{_generic_manifolds_Definition_} 
Let $M$ be a complex manifold admitting a hyperk\"ahler
structure $\c H$. We say that $M$ is {\bf of general type}
or {\bf generic} with respect to $\c H$ if all elements of the group
\[ \bigoplus\limits_p H^{p,p}(M)\cap H^{2p}(M,\Z)\subset H^*(M)\] 
are $SU(2)$-invariant. 
We say that $M$ is {\bf Mumford--Tate generic} 
if for all $n\in \Z^{>0}$, all the cohomology classes
\[ \alpha \in
   \bigoplus\limits_p H^{p,p}(M^n)\cap H^{2p}(M^n,\Z)\subset H^*(M^n)
\] 
are $SU(2)$-invariant. In other words,
$M$ is Mumford--Tate generic if for all
$n\in {\Bbb Z}^{>0}$, the $n$-th power $M^n$ is
generic. Clearly, Mumford--Tate generic
implies generic.

\hfill

\proposition \label{_generic_are_dense_Proposition_} 
Let $M$ be a compact manifold, $\c H$ a hyperk\"ahler
structure on $M$ and $S$
be the set of induced complex structures over $M$. Denote by 
$S_0\subset S$ the set of $L\in S$ such that 
$(M,L)$ is Mumford-Tate generic with respect to $\c H$. 
Then $S_0$ is dense in $S$. Moreover, the complement
$S\backslash S_0$ is countable.

{\bf Proof:} This is Proposition 2.2 from
\cite{_Verbitsky:Symplectic_II_}
\endproof

\hfill

\ref{_G_M_invariant_implies_trianalytic_Theorem_} has the following
immediate corollary:

\corollary \label{_hyperkae_embeddings_Corollary_} 
Let $M$ be a compact holomorphically symplectic 
manifold. Assume that $M$ is of general type with respect to
a hyperk\"ahler structure $\c H$.
Let $S\subset M$ be closed complex analytic
subvariety. Then $S$ is trianalytic 
with respect to $\c H$.

\endproof

\hfill

In \cite{_Verbitsky:hypercomple_},
\cite{_Verbitsky:Desingu_}, \cite{_Verbitsky:DesinguII_},
we gave a number of equivalent definitions of a singular hyperk\"ahler
and hypercomplex variety. We refer the reader to 
\cite{_Verbitsky:DesinguII_} for the precise definition;
for our present purposes it suffices to say that all
trianalytic subvarieties are hyperk\"ahler varieties.
The following Desingularization Theorem is very 
useful in the study of trianalytic subvarieties.

\hfill

\theorem \label{_desingu_Theorem_}
(\cite{_Verbitsky:DesinguII_})
Let $M$ be a hyperk\"ahler or a hypercomplex variety,
$I$ an induced complex structure.
Let \[ \widetilde{(M, I)}\stackrel n\arrow (M,I)\] 
be the normalization of
$(M,I)$. Then $\widetilde{(M, I)}$ is smooth and
has a natural hyperk\"ahler
structure $\c H$, such that the associated
map $n:\; \widetilde{(M, I)} \arrow (M,I)$ agrees with $\c H$.
Moreover, the hyperk\"ahler 
manifold $\tilde M:= \widetilde{(M, I)}$
is independent from the choice of induced complex structure $I$.

\endproof

\subsection{Simple hyperk\"ahler manifolds}

\definition \label{_simple_hyperkahler_mfolds_Definition_} 
(\cite{_Beauville_}) A connected simply connected 
compact hy\-per\-k\"ah\-ler manifold 
$M$ is called {\bf simple} if $M$ cannot be represented 
as a product of two hyperk\"ahler manifolds:
\[ 
   M\neq M_1\times M_2,\ \text{where} \ dim\; M_1>0 \ dim\; M_2>0
\] 
Bogomolov proved that every compact hyperk\"ahler manifold has a finite
covering which is a product of a compact torus
and several simple hyperk\"ahler manifolds. 
Bogomolov's theorem implies the following result
(\cite{_Beauville_}):

\hfill

\theorem\label{_simple_mani_crite_Theorem_}
Let $M$ be a compact hyperk\"ahler manifold.
Then the following conditions are equivalent.
\begin{description}
\item[(i)] $M$ is simple
\item[(ii)] $M$ satisfies $H^1(M, \R) =0$, $H^{2,0}(M) =\C$,
where $H^{2,0}(M)$ is the space of $(2,0)$-classes taken with
respect to any of induced complex structures.
\end{description}

\endproof


\subsection{Hyperholomorphic bundles}
\label{_hyperholo_Subsection_}


This subsection contains several versions of a
definition of hyperholomorphic connection in a complex
vector bundle over a hyperk\"ahler manifold.
We follow \cite{_Verbitsky:Hyperholo_bundles_}.

\hfill

 Let $B$ be a holomorphic vector bundle over a complex
manifold $M$, $\nabla$ a  connection 
in $B$ and $\Theta\in\Lambda^2\otimes End(B)$ be its curvature. 
This connection
is called {\bf compatible with a holomorphic structure} if
$\nabla_X(\zeta)=0$ for any holomorphic section $\zeta$ and
any antiholomorphic tangent vector $X$. If there exist
a holomorphic structure compatible with the given
Hermitian connection then this connection is called
{\bf integrable}.
 
\hfill
 
One can define a {\bf Hodge decomposition} in the space of differential
forms with coefficients in any complex bundle, in particular,
$End(B)$.

\hfill

\theorem \label{_Newle_Nie_for_bu_Theorem_}
Let $\nabla$ be a Hermitian connection in a complex vector
bundle $B$ over a complex manifold. Then $\nabla$ is integrable
if and only if $\Theta\in\Lambda^{1,1}(M, \End(B))$, where
$\Lambda^{1,1}(M, \End(B))$ denotes the forms of Hodge
type (1,1). Also,
the holomorphic structure compatible with $\nabla$ is unique.

{\bf Proof:} This is Proposition 4.17 of \cite{_Kobayashi_}, 
Chapter I.
$\blacksquare$

\hfill

\definition \label{_hyperho_conne_Definition_}
Let $B$ be a Hermitian vector bundle with
a connection $\nabla$ over a hyperk\"ahler manifold
$M$. Then $\nabla$ is called {\bf hyperholomorphic} if 
$\nabla$ is
integrable with respect to each of the complex structures induced
by the hyperk\"ahler structure. 
 
As follows from 
\ref{_Newle_Nie_for_bu_Theorem_}, $\nabla$ is hyperholomorphic
if and only if its curvature $\Theta$ is of Hodge type (1,1) with
respect to any of complex structures induced by a hyperk\"ahler 
structure.

As follows from \ref{_SU(2)_inva_type_p,p_Lemma_}, 
$\nabla$ is hyperholomorphic
if and only if $\Theta$ is a $SU(2)$-invariant differential form.

\hfill

\example \label{_tangent_hyperholo_Example_} 
(Examples of hyperholomorphic bundles)

\begin{description}

\item[(i)]
Let $M$ be a hyperk\"ahler manifold, $TM\otimes_\R \C$ 
a complexification of its tangent bundle
equipped with Levi--Civita connection $\nabla$. Then $\nabla$
is integrable with respect to each induced complex structure,
and hence, Yang--Mills.

\item[(ii)] For $B$ a hyperholomorphic bundle, all its tensor powers
are also hyperholomorphic.

\item[(iii)] Thus, the bundles of differential forms on a hyperk\"ahler
manifold are also hyperholomorphic.

\end{description}


\subsection{Stable bundles and Yang--Mills connections.}


This subsection is a compendium of the most
basic results and definitions from the Yang--Mills theory
over K\"ahler manifolds, concluding in the fundamental
theorem of Uhlenbeck--Yau \cite{_Uhle_Yau_}.

\hfill

\definition\label{_degree,slope_destabilising_Definition_} 
Let $F$ be a coherent sheaf over
an $n$-dimensional compact K\"ahler manifold $M$. We define
$\deg(F)$ as
 
\[ 
   \deg(F)=\int_M\frac{ c_1(F)\wedge\omega^{n-1}}{vol(M)}
\] 
and $\text{slope}(F)$ as
\[ 
   \text{slope}(F)=\frac{1}{\text{rank}(F)}\cdot \deg(F). 
\]
The number $\text{slope}(F)$ depends only on a
cohomology class of $c_1(F)$. 

Let $F$ be a coherent sheaf on $M$
and $F'\subset F$ its proper subsheaf. Then $F'$ is 
called {\bf destabilizing subsheaf} 
if $\text{slope}(F') \geq \text{slope}(F)$

A holomorphic vector bundle $B$ is called {\bf stable}
\footnote{In the sense of Mumford-Takemoto}
if it has no destabilizing subsheaves.
 
\hfill

Later on, we usually consider the bundles $B$ with $deg(B)=0$.

\hfill

Let $M$ be a K\"ahler manifold with a K\"ahler form $\omega$.
For differential forms with coefficients in any vector bundle
there is a Hodge operator $L: \eta\arrow\omega\wedge\eta$.
There is also a fiberwise-adjoint Hodge operator $\Lambda$
(see \cite{_Griffi_Harri_}).
 
\hfill

\definition \label{Yang-Mills_Definition_}
Let $B$ be a holomorphic bundle over a K\"ahler manifold $M$
with a holomorphic Hermitian connection $\nabla$ and a 
curvature $\Theta\in\Lambda^{1,1}\otimes End(B)$.
The Hermitian metric on $B$ and the connection $\nabla$
defined by this metric are called {\bf Yang-Mills} if 

\[
   \Lambda(\Theta)=constant\cdot \Id\restrict{B},
\]
where $\Lambda$ is a Hodge operator and $\Id\restrict{B}$ is 
the identity endomorphism which is a section of $End(B)$.

Further on, we consider only these Yang--Mills connections
for which this constant is zero.

\hfill

A holomorphic bundle is called  {\bf indecomposable} 
if it cannot be decomposed onto a direct sum
of two or more holomorphic bundles.

\hfill

The following fundamental 
theorem provides examples of Yang-\--Mills \linebreak bundles.

\theorem \label{_UY_Theorem_} 
(Uhlenbeck-Yau)
Let B be an indecomposable
holomorphic bundle over a compact K\"ahler manifold. Then $B$ admits
a Hermitian Yang-Mills connection if and only if it is stable, and
this connection is unique.
 
{\bf Proof:} \cite{_Uhle_Yau_}. \endproof

\hfill

\proposition \label{_hyperholo_Yang--Mills_Proposition_}
Let $M$ be a hyperk\"ahler manifold, $L$
an induced complex structure and $B$ be a complex vector
bundle over $(M,L)$. 
Then every hyperholomorphic connection $\nabla$ in $B$
is Yang-Mills and satisfies $\Lambda(\Theta)=0$
where $\Theta$ is a curvature of $\nabla$.
 
\hfill

{\bf Proof:} We use the definition of a hyperholomorphic 
connection as one with $SU(2)$-invariant curvature. 
Then \ref{_hyperholo_Yang--Mills_Proposition_}
follows from the

\hfill

\lemma \label{_Lambda_of_inva_forms_zero_Lemma_}
Let $\Theta\in \Lambda^2(M)$ be a $SU(2)$-invariant 
differential 2-form on $M$. Then
$\Lambda_L(\Theta)=0$ for each induced complex structure
$L$.\footnote{By $\Lambda_L$ we understand the Hodge operator 
$\Lambda$ associated with the K\"ahler complex structure $L$.}

{\bf Proof:} This is Lemma 2.1 of \cite{_Verbitsky:Hyperholo_bundles_}.
\endproof
 
\hfill

Let $M$ be a compact hyperk\"ahler manifold, $I$ an induced 
complex structure. 
For any stable holomorphic bundle on $(M, I)$ there exists a unique
Hermitian Yang-Mills connection which, for some bundles,
turns out to be hyperholomorphic. It is possible to tell when
this happens (though in the present paper we never use
this knowledge).

\hfill

\theorem
Let $B$ be a stable holomorphic bundle over
$(M,I)$, where $M$ is a hyperk\"ahler manifold and $I$
is an induced complex structure over $M$. Then 
$B$ admits a compatible hyperholomorphic connection if and only
if the first two Chern classes $c_1(B)$ and $c_2(B)$ are 
$SU(2)$-invariant.\footnote{We use \ref{_SU(2)_commu_Laplace_Lemma_}
to speak of action of $SU(2)$ in cohomology of $M$.}

{\bf Proof:} This is Theorem 2.5 of
 \cite{_Verbitsky:Hyperholo_bundles_}. \endproof


\section{Trianalytic subvarieties of 
powers of K3 surfaces}
\label{_appendix_subva_produ_Section_}


\subsection{Trianalytic subvarieties of a product of a K3 surface with itself}

Let $M$ be any manifold, $M^n = M\times ... \times M$ its $n$-th
product with itself. We define the ``natural'' subvarieties
of $M$, recursively, as follows.

\begin{equation}\label{_recu_subva_in_M^n_Equation_}
\begin{minipage}[m]{0.8\linewidth}
\begin{description}
\item[(i)] Natural subvarieties of $M$ are $M$ and points.
\item[(ii)] Let $Z\subset M^n$ by a natural subvariety.
The following subvarieties of $M^{n+1}$ are natural.
\begin{description}
\item[a.] $Z_M := Z\times M \subset M^n \times M = M^{n+1}$
\item[b.] $Z_t:= Z\times \{t\} \subset M^n \times M= M^{n+1}$,
depending on a point $t\in M$.
\item[c.] $Z_i := 
  \bigg\{ (m_1, ... m_{n+1}) \in Z\times M \;\;|\;\; 
    m_i=m_{n+1} \bigg\}$ depending on a number
    $i\in \{ 1, ...,  n\}$
\end{description}
\end{description}
\end{minipage}
\end{equation}

The main result of this section is the following theorem.

\hfill

\theorem \label{_subva_product_Therem_}
Let $M$ be a hyperk\"ahler K3-surface which has no 
hyperk\"ahler automorphisms, and $X\subset M^{n}$ an
irreducible trianalytic subvariety.
Then $X$ is ``natural'', in the sense of 
\eqref{_recu_subva_in_M^n_Equation_}.

\hfill

{\bf Proof:} Let \[ \check\Pi_{n+1}:\; M^{n+1}\arrow M^n\]
be the natural projection $m_1, ..., m_{n+1} \arrow
m_1, ..., m_n$. Clearly,
$\check\Pi_{n+1}(X)$ is irreducible and trianalytic.
Using induction, we may assume that 

\begin{equation}\label{_indu_assu_for_subva_in_product_Equation_}
\begin{minipage}[m]{0.8\linewidth} a trianalyiic subvariety
$X\subset M^k$ is natural, in the sense 
of \eqref{_recu_subva_in_M^n_Equation_},
for $k\leq n$.
\end{minipage}
\end{equation}
Clearly, by \eqref{_indu_assu_for_subva_in_product_Equation_}
$\check\Pi_{n+1}(X)$ is 
of the type \eqref{_recu_subva_in_M^n_Equation_}.
All varieties $Z$ of type \eqref{_recu_subva_in_M^n_Equation_}
are isomorphic to $M^k$, for $k=\dim_{\Bbb H} Z$. Thus,
$X$ is realized as a subvariety in 
\[ \check\Pi_{n+1}(X)\times M =
   M^{\dim_{\Bbb H}\check\Pi_{n+1}(X)+1}.
\]
Unless $\dim_{\Bbb H} \check\Pi_{n+1}(X) =n$,
\eqref{_indu_assu_for_subva_in_product_Equation_} implies that
$X$ is a ``natural'' subvariety.
Thus, to prove \eqref{_subva_product_Therem_},
we may assume that $\check\Pi_{n+1}(X) = M^n$.

For a point $t\in M$, let
$X_t:= \bigg\{ (m_1, ... m_{n+1}) \in X \;\; |\;\; m_{n+1}=t \bigg\}$.
The subvariety $X_t\subset M^n$ is not necessarily irreducible,
because the components of $X_t$ may ``flow together''
as $t$ changes, so that $X = \bigcup_{t\in M} X_t$ is still
irreducible, while the $X_t$'s are not.
However, all components of $X_t$ must be
deformationally equivalent in the class of ``natural'' subvarieties
of $M^n$, in order for $X$ to be irreducible.

Since $X$ is irreducible and $X= \cup_t X_t$, 
then either $X_t = X$ for one value of $t$ (and in this case
$X= M^{n}\times \{t\}$),
or $X_t \neq\emptyset$ for $t$ in a positive-dimensional trianalytic
subset of $M$. Since $\dim_{\Bbb H} M =1$, this subset coinsides
with $M$. Using \eqref{_indu_assu_for_subva_in_product_Equation_}, 
we obtain that all irreducible components of $X_t$ are of the type
\eqref{_recu_subva_in_M^n_Equation_}. All ``natural''
subvarieties of $M^n$ of complex codimension $2=\dim_\C M$ are
given by either $m_i=m_j$ for some distinct fixed indices $i, j$,
or by $m_i=t$ for a fixed index $i$ and a fixed point
$t\in M$. We proceed on case-by-case basis.

\begin{description}
\item[(i)] For some $t$, $X_t$ contains a component
$X_t^{i,j}$ given by $m_i=m_j$ for  distinct fixed indices $i, j$.
Since $X_t^{i,j}$ is rigid in the class of natural subvarieties,
and $X$ is irreducible, this implies that
$X_t$ contains $M^n_{i,j}\times\{t\}$ for all $t$, where
$M^n_{i,j}\subset M^n$ is a subvariety given by $m_i=m_j$. Since
$\dim X= M^n_{i,j} \times M$, $X$ irreducible
and $M^n_{i,j} \times M\subset X$, this implies that
$X=M^n_{i,j} \times M$. This proves 
\ref{_subva_product_Therem_} (case (i)).

\item[(ii)] For some $t$, $X_t$ contains a component
$X^i_t(m)$, given by $m_i=m$, for a fixed index $i$ and a fixed point
$m\in M$. Deforming $X^i_t(m)$ in the class of natural subvarieties,
we obtain again $X^i_t(m')$, with different $m'$. Taking 
a union of all  $X^i_t(m)\subset X_t$, for some fixed $i$
and varying $t$ and $m$, we obtain a closed subvariety of $X$
of the same dimension as $X$. Since $X$ is irreducible,
{\it all} components of $X_t$ are given by 
$m_i=m$, for a fixed index $i$ and a fixed point
$m\in M$. Consider a trianalytic subvariety $Z\subset M^2$,
\[ 
  Z:= \{ (m, t) \in M^2 \;\;|\;\; X^i_t(m)\subset X_t \} 
\]
To prove \ref{_subva_product_Therem_} (case (ii)),
it suffices to show that $Z$ is natural, in the sense of
\eqref{_recu_subva_in_M^n_Equation_}. We reduced
\ref{_subva_product_Therem_} to the case of trianalytic
subvarieties in $M^2$.
\end{description}

The following lemma finishes the proof of \ref{_subva_product_Therem_}.

\hfill

\lemma \label{_subva_M^2_Lemma_}
Let $M$ be a hyperk\"ahler K3-surface which has no 
hyperk\"ahler automorphisms, and
 $X\subset M^2$ a closed irreducible trianalytic
subvariety of $M^2$. Then $X$ is a ``natural''
subvariety of $M^2$, in the sense of 
\eqref{_recu_subva_in_M^n_Equation_}.

\hfill

{\bf Proof:} Let $\pi_1, \pi_2:\; M^2 \arrow M$ be the 
natural projections. Assume that neither $\pi_1(X)$
nor $\pi_2(X)$ is a point, and $X\subsetneqq M^2$ 
(otherwise $X$ obviously satisfies 
\eqref{_recu_subva_in_M^n_Equation_}). 
Let $\tilde X\stackrel n
\arrow X$ be the desingularization of $X$,
given by \ref{_desingu_Theorem_}. 

Consider the maps $p_i:\; \tilde X \arrow M$,
$i= 1,2$, given by $p_i:= n\circ \pi_i$. Since
$\dim X = \dim M$, and $p_i$ is non-trivial,
these maps have non-degenerate Jacobians in general point.
Fix an induced complex structure $I$ on $M$, and consider
$X$, $\tilde X$, $M^2$ as complex varieties and
$p_i$ as holomorphic maps. Let $\Theta_M$ be the holomorphic 
symplectic form of $M$. Then $p_i^* \Theta_M$ gives a section
of the canonical class of $\tilde X$. Since $\tilde X$ is 
compact and hyperk\"ahler, any non-zero section of the canonical class
is nowhere degenerate. Thus, $p_i^* \Theta_M$ is nowhere degenerate,
and the Jacobian of $p_i$ nowhere vanishes. Therefore,
$p_i$ is a covering. Since $X$ is irreducible, $\tilde X$
is connected, and since $M$ is simply connected, 
$p_i$ is isomorphism. Since $M$ has no hyperk\"ahler
automorphisms, except identity, $X$ is a graph
of an identity map. 
This proves \ref{_subva_M^2_Lemma_} and 
\ref{_subva_product_Therem_}.
\endproof

\hfill

\corollary \label{_complex_subva_M^n_Corollary_}
Let $M$ be a complex K3 surface with no complex
automorphisms. Assume that $M$ admits a hyperk\"ahler
structure $\c H$ such that $M$ is Mumford-Tate generic
with respect to $\c H$. Let $X$ be an irreducible
complex subvariety of $M^n$. Then $X$ is a ``natural''
subvariety of $M^n$, in the sense of 
\eqref{_recu_subva_in_M^n_Equation_}.

{\bf Proof:} 
By \ref{_hyperkae_embeddings_Corollary_}, 
$X$ is trianalytic. Now \ref{_complex_subva_M^n_Corollary_} is implied by
\ref{_subva_product_Therem_}. \endproof

\subsection{Subvarieties of symmetric powers of varieties}
\label{_subva_of_symme_special_Subsection_}

In this section, we fix the notation regarding
the ``natural'' subvarieties of the symmetric powers of 
complex varieties. 

\hfill

Let $M$ be a complex variety and 
$M^{(n)}$ its symmetric power,
$M^{(n)}=M^n/S_n$. The space $M^{(n)}$ has 
a natural stratification by {\bf diagonals} $\Delta_{(\alpha)}$,
which 
are numbered by Young diagrams 
\[ \alpha = (n_1 \geq n_2 ...\geq n_k), \ \ \sum n_i = n.\] 
This stratification is constructed as follows.
Let $\sigma:\; M^n \arrow M^{(n)}$ be the natural finite
map (a quotient by the symmetric group). Then
\begin{equation} \label{_Del-a_alpha_definition_Equation_}
\begin{split}
  \Delta_{(\alpha)} := &
  \sigma \left (\left\{ (x_1, x_2, ... , x_n)\in M^n \ \ \bigg | \;\; 
\right.\right.\\ &
              x_1 = x_2 = ... = x_{n_1}, \ \  
 	       x_{n_1+1} = x_{n_1+2} = ... = x_{n_1+n_2}, \  
 	       \\ &
         \left.\left. \vphantom{\ \ \bigg | \;\;}
              ..., \ \  x_{\sum_{i=1}^{k-1} n_i+1} = 
	       x_{\sum_{i=1}^{k-1} n_i+2} = ... = x_n \right\} \right)
\end{split}
\end{equation}
where $\sigma:\; M^n \arrow M^{(n)}$ is the natural quotient map.

\hfill

 Consider
a Young diagram $\alpha$,
\[ 
   \alpha = (n_1 \geq n_2 ...\geq n_k), \ \ \sum n_i = n.
\] 
As in \eqref{_Del-a_alpha_definition_Equation_}, 
$\alpha$ corresponds to a diagonal $\Delta_{(\alpha)}$, which is a
closed subvariety of $M^{(n)}$.
Fix a subset $\c A \subset \{1, ..., k\}$, and let
$\phi:\; \c A \arrow M$ be an arbitrary 
map.  Then $\Delta_{(\alpha)}(\c A, \phi)\subset \Delta_{(\alpha)}$
is given by

\begin{equation} \label{_Delta(A,phi)_alpha_definition_Equation_}
\begin{split}
  \Delta_{(\alpha)}(\c A, \phi) := &
  \sigma \left (\left\{ (x_1, x_2, ... , x_n)\in M^n \ \ \bigg | \;\; 
\right.\right.\\ &
              x_1 = x_2 = ... = x_{n_1}, \ \  
               x_{n_1+1} = x_{n_1+2} = ... = x_{n_1+n_2}, \  
               \\ &
              ..., \ \  x_{\sum_{i=1}^{k-1} n_i+1} = 
               x_{\sum_{i=1}^{k-1} n_i+2} = ... = x_n,\\
        & \text{\ \ and\ \ } \forall i \in \c A, \ x_i = \phi(i)
\left.\left. \vphantom{\ \ \bigg | \;\;} \right\} \right),
\end{split}
\end{equation}
where $\sigma:\; M^n \arrow M^{(n)}$ is the standard 
projection. For ``sufficiently generic'' K3 surfaces,
all complex subvarieties in $M^n$ are given by
\ref{_complex_subva_M^n_Corollary_}. From 
\ref{_complex_subva_M^n_Corollary_}, it is easy to deduce
the following result.

\hfill

\proposition\label{_subva_in_M^(n)_Proposition_}
Let $M$ be a complex K3 surface with no complex
automorphisms. Assume that $M$ admits a hyperk\"ahler
structure $\c H$ such that $M$ is Mumford-Tate generic
with respect to $\c H$. Let $X$ be an irreducible
complex subvariety of $M^n$.
Then $X= \Delta_{(\alpha)}(\c A, \phi)$
for appropriate $\alpha, \c A, \phi$.

\endproof


\section{Hilbert scheme of points}
\label{_Hilbert_sche_Section_}


For the definitions and results related to the 
 Hilbert scheme of points on a surface, see the excellent
lecture notes of H. Nakajima \cite{_Nakajima_}. 

\subsection{Symplectic structure of the Hilbert scheme}
\label{_symple_on_Hi_Subsection_}

\definition
Let $M$ be a complex surface. The $n$-th
{\bf Hilbert scheme of points}, also called {\bf Hilbert scheme}
of $M$ (denoted by $M^{[n]})$ is the scheme classifying the 
zero-dimensional subschemes of $M$ of length $n$.

\hfill

There is a natural projection 
$\pi:\; M^{[n]} \arrow M^{(n)}$ associating to a 
subscheme $S\subset M$ the set of points $x_i\in Sup(S)$
of support of $S$, taken with the multiplicities equal to the 
length of $S$ in $x_i$.

It is well-known that
a Hilbert scheme of a smooth surface is a smooth manifold, 
and the fibers of $\pi$ are irreducible and reduced
(see, e. g. \cite{_Nakajima_}).

\hfill

For our purposes, the most important property of the Hilbert
scheme is the existence of the non-degenerate holomorphic symplectic
structure, for $M$ holomorphically symplectic.

\hfill

Let $X$ be an irreducible 
complex analytic space, which is reduced in generic point, and
$\Omega^1 X$ the sheaf
of Ka\"hler differentials on $X$. We denote the exterior square
$\Lambda^2_{\calo_X}\Omega^1 X$ by $\Omega^2 X$. The sections
of $\Lambda^2_{\calo_X}\Omega^1 X$ are called {\bf 2-forms
on $X$}. 
We say that two-forms $\omega_1, \omega_2$ are {\bf
equal up to a torsion} if $\omega_1=\omega_2$ in the generic
point of $X$.

\hfill

\proposition \label{_Hilbert_symple_Proposition_}
(Beauville)
Let $M$ be a smooth complex surface equipped with a nowhere
degenerate holomorphic 2-form. Then 
\begin{description}
\item[(i)] the Hilbert scheme $M^{[n]}$
is equipped with a natural, nowhere degenerate holomorphic 2-form
$\Theta_{M^{[n]}}$. 
\item[(ii)] Consider the Cartesian square
\begin{equation}\label{_M^[n]_Cartesian_Equation_}
\begin{CD}
\tilde M^{[n]} @>{\tilde \pi}>>  M^{n} \\
@VV{\tilde \sigma}V @VV{\sigma}V \\
M^{[n]} @>{\pi}>>  M^{(n)}
\end{CD}
\end{equation}
Let $\Theta_{M^{n}}$ be the natural symplectic form on $M^n$.
Then the pullback $\tilde\sigma^* \Theta_{M^{[n]}}$
is equal to the pullback $\tilde \pi^* \Theta_{M^{n}}$, outside
of the subvariety $D\subset \tilde M^{[n]}$ of codimension 2.
\item[(iii)] The complex analytic space $\tilde M^{[n]}$
is irreducible, and $\tilde\sigma^* \Theta_{M^{[n]}}$
is  equal up to a torsion to $\tilde \pi^* \Theta_{M^{n}}$.
\end{description}

\hfill

{\bf Proof:} In \cite{_Beauville_}, A. Beauville proved 
the conditions (i) and (ii).  Clearly, (ii) implies that
$\tilde\sigma^* \Theta_{M^{[n]}}$
is  equal up to a torsion to $\tilde \pi^* \Theta_{M^{n}}$,
assuming that $\tilde M^{[n]}$ is irreducible. It remains
to show that $\tilde M^{[n]}$ is irreducible. 

\hfill

The following argument can be easily generalized to
a more general type of Cartesian squares. We only use that
the arrow $\pi$ is rational, $\sigma$ is finite
and generically etale, and the
varieties $M^n$, $M^{[n]}$ are irreducible.

\hfill

Let $U$ be the general open stratum of $\tilde M^{[n]}$,
\[ U := 
\tilde M^{[n]}\backslash \bigcup_\alpha \tilde \sigma^{-1}\Delta_{[\alpha]}
\]
The map $\tilde \pi:\; U \arrow M^n$ is an open embedding.
Therefore, the variety $U$ is irreducible.
To prove that $\tilde M^{[n]}$ is irreducible, we need to show
that for all points $x\in \tilde M^{[n]}$, there exists a sequence
$\{x_i\}\subset U$ which converges to $x$. 
Since $M^{[n]}$ is irreducible, there exists a sequence
$\{ \underline{x_i}\}\in \tilde\sigma(U)$ converging to $\tilde \sigma(x)$.
Consider the sequence $\{\pi(\underline{x_i})\} \subset M^{(n)}$. 
The general stratum $\tilde \pi(U)$ of $M^n$ is identified with $U$,
since $\tilde \pi\restrict U$ is an isomorphism.
Lifting $\{\pi(\underline{x_i})\}$ to $M^n$,
we obtain a sequence $\{x_i\}\subset \tilde \pi(U)=U$.
Taking a subsequence of $\{x_i\}$, we can assure that
it converges to a point in a finite set 
$\tilde\sigma^{-1}(\tilde\sigma(x))$. 
By an appropriate choice of the lifting,
we obtain a sequence converging to any point
in $\tilde\sigma^{-1}(\tilde\sigma(x))$,
in particular, $x$. This proves that 
$\tilde M^{[n]}$ is irreducible.
 \endproof

\hfill

\remark
{}From \ref{_Hilbert_symple_Proposition_} and Calabi-Yau theorem
(\ref{_symplectic_=>_hyperkahler_Proposition_}), it follows
immediately that $M^{[n]}$ admits a hyperk\"ahler structure,
if $M$ is a K3 surface or a compact torus.
However, this hyperk\"ahler structure is not in any
way related to the hyperk\"ahler structures on $M$.

\hfill

\remark 
In the preliminary version of \cite{_Nakajima_}, it was stated without proof 
that  $\tilde\sigma^* \Theta_{M^{[n]}} =\tilde \pi^* \Theta_{M^{n}}$.
This statement seems to be subtle, and I was unable to find the proof.
However, a weaker version of this equality can be proven.

\hfill

\proposition \label{_simplec_on_subva_Proposition_}
Let $A_{M^{[n]}}:= 
\tilde\sigma^* \Theta_{M^{[n]}}$, $A_{M^n}:= \tilde \pi^* \Theta_{M^{n}}$,
be the forms defined in \ref{_Hilbert_symple_Proposition_}.
Then for all closed 
complex subvarieties $X\stackrel i \hookrightarrow \tilde M^{[n]}$,
the 2-forms $i^* A_{M^{[n]}}$, $i^*A_{M^n}\in \Omega^2 X$ 
are equal outside of singularities of $X$.

\hfill

{\bf Proof:}\footnote{The proof is based on ideas of D. Kaledin.}
The forms  $A_{M^{[n]}}$, $A_{M^n}$ are equal up to torsion.
Therefore, their difference lies in the torsion subsheaf
of $\Omega^2 \tilde M^{[n]}$.
To prove that the 2-forms $i^* A_{M^{[n]}}$, $i^*A_{M^n}\in \Omega^2 X$ 
are equal outside of singularities of $X$ it suffices to show the
following: for all torsion sections $\omega \in \Omega^2 \tilde M^{[n]}$,
the pullback $i^*\omega$ lies in the torsion of $\Omega^2 X$ 

Let $\Delta_{[\alpha]}$ be a stratum of $M^{[n]}$,
defined by a Young diagram $\alpha$ as in 
Subsection \ref{_Young_Hilb_Subsection_}, and 
$X\stackrel i \hookrightarrow \tilde M^{[n]}$ 
an irreducible component of $\tilde \sigma^{-1}(X)$, considered as
a complex subvariety of $\tilde M^{[n]}$. 
As a first step in proving \ref{_simplec_on_subva_Proposition_},
we show that for all torsion sections $\omega \in \Omega^2 \tilde M^{[n]}$,
the pullback $i^*\omega$ lies in the torsion of $\Omega^2 X$, for 
this particular choice of $X$. Consider a generic point
 $x\in\Delta_{[\alpha]}$, and let
$V\subset M^{[n]}$ be a neighbourhood of $x$ in
$M^{[n]}$,  $U\subset \Delta_{[\alpha]}$ be a neighbourhood
of $x\in\Delta_{[\alpha]}$.
For an appropriate choice of $V$, $U$, these varieties
are equipped with a locally trivial fibration $V \stackrel p\arrow U$,
 inverse to a natural embedding $U\hookrightarrow V$.
Assume also that $U$ consists entirely of generic points
of $\Delta_{[\alpha]}$.
Let $\tilde x$ be a point of $\tilde \sigma^{-1}(x) \cap X$,
and $\tilde V$ be a component
of $\tilde \sigma^{-1}(V)$ which contains $\tilde x$,
and $\tilde U:= \tilde \sigma^{-1}(U) \cap \tilde V$.
Since $U$ consists of generic points of $\Delta_{[\alpha]}$,
and $\tilde \sigma$ is finite,
the map $\tilde \sigma:\; \tilde U \arrow U$ 
is etale. Therefore, $\tilde V$ is equipped
with a locally trivial 
fibration $\tilde V \stackrel {\tilde p}\arrow \tilde U$,
inverse to a natural embedding $\tilde U\hookrightarrow\tilde  V$.
Using this fibration, we decompose the 
sheaf of differentials on $\tilde V$ as follows:
\[
    \Omega^1 \tilde V = {\tilde p}^*\Omega^1 \tilde U 
    \oplus \Omega^1_{\tilde p} \tilde V
\]
where $\Omega^1_{\tilde p} \tilde V$ is the sheaf of relative differentials
of $\tilde V$ along $\tilde p$. Clearly, for all sections  
$\omega \in \Omega^1_{\tilde p} \tilde V$,
 the pullback of $\omega$ under $\tilde U \hookrightarrow\tilde  V$
is zero. On the other hand, $\tilde U$ is smooth, and therefore
${\tilde p}^*\Omega^1 \tilde U$ has no torsion. Thus,
the torsion-component of $\Omega^1 \tilde V$ is contained in
$\Omega^1_{\tilde p} \tilde V$ and vanishes on
$\tilde U$. A similar argument implies that
the torsion-component of $\Omega^2 \tilde V$ is contained in
\[
  \Omega^2_{\tilde p} \tilde V\oplus 
  {\tilde p}^*\Omega^1 \tilde U \otimes \Omega^1_{\tilde p} \tilde V
  \subset \Omega^2 \tilde V
\]
and also vanishes on $\tilde U$. Therefore, all torsion 
components on $\Omega^2 \tilde M^{[n]}$ vanish on
$X$, where $X$ is an irreducible component of the preimage of the
stratum of $M^{[n]}$. Consider a stratification
of $\tilde M^{[n]}$ by such $X$'s. For any
subvariety $Y\subset \tilde M^{[n]}$, let
$\tilde \Delta_{[\alpha]}$ be the smallest stratum
of $\tilde M^{[n]}$ which contains $Y$. Then, the
set $Y_g$ of generic points of $Y$ is contained in the set
$\tilde U_{[\alpha]}$ of the generic points of 
$\tilde \Delta_{[\alpha]}$. But, as we have seen,
the restrictions of the forms $A_{M^{[n]}}$, $A_{M^n}$ to
$\tilde U_{[\alpha]}$ coinside. Therefore,
restrictions of $A_{M^{[n]}}$, $A_{M^n}$ to 
$Y_g \subset \tilde U_{[\alpha]}$ are equal.
This proves \ref{_simplec_on_subva_Proposition_}.
\endproof

\hfill

In the situation similar to the above,
we will say that the forms $A_{M^{[n]}}$, $A_{M^n}$
are {\bf equal on subvarieties}. 

\hfill

From the fact that $A_{M^{[n]}}$, $A_{M^n}$
are equal on subvarieties we immediately obtain the following.

\hfill

\claim \label{_triana_finite_in_gene_Claim_}
Let $X\subset M^{[n]}$ be a complex subvariety of the
Hilbert scheme. Assume that the holomorphic symplectic
form is non-degenerate in the generic point of $X$
(this happens, for instance, when $X$ is trianalytic).
Then the restriction $\pi\restrict X$ of 
$\pi:\;M^{[n]}\arrow M^{(n)}$ to $X$ is 
finite in generic point of $X$

{\bf Proof:} This result is easily implied by 
\ref{_simplec_on_subva_Proposition_}. 
For details of the proof, the 
reader is referred to \ref{_holo_symple_semismall_Proposition_},
\ref{_triana_speci_semismall_Corollary_}.
\endproof

\hfill

The rest of this section is not used directly anywhere in this
paper. A reader who does not like perverse sheaves is invited to 
skip the rest and proceed to Section \ref{_unive_subva_Section_}.

\subsection{Cohomology of the Hilbert scheme}

For perverse sheaves, we freely use terminology and results of 
\cite{_Asterisque_100_} . For the computation of cohomology of the Hilbert
scheme via perverse sheaves, see \cite{_Nakajima_}.

\hfill

\definition
Let $X$ be an irreducible complex variety and $F$ a perverse sheaf on
$M$. The $F$ is called {\bf a Goresky-MacPherson sheaf},
or {\bf a sheaf of GM-type} if
it has no proper subquotient perverse sheaves
with support in $Z\subsetneq X$. For an arbitrary 
perverse sheaf $F$, consider the Goresky-MacPherson
subquotient $F_{GM}$ of $F$, such that for 
a nonempty Zariski open set $U\subset X$, 
$F\restrict U = F_{GM}\restrict U$.
Such subquotient is obviously unique; we 
call it the Goresky-MacPherson extension of $F$.
The Intersection Cohomology sheaf $IC(X)$ is the
Goresky-MacPherson extension of the constant sheaf
$\C_X$.

\hfill

\definition
Let $X$ be a complex variety. The $X$ is called
{\bf homology rational} if the constant sheaf $\C_X$ on $X$,
considered as a complex of sheaves, is isomorphic to
the Intersection Cohomology perverse sheaf. 
The variety $X$ is called {\bf weakly homology rational}
if the Intersection Cohomology sheaf (considered 
as a complex of sheaves) is constructible (i. e.,
the cohomology of this complex of sheaves
are zero in all but one degree).

\hfill

\remark
Clearly, for a homology rational variety,
the Intersection Cohomology coinsides with the
standard cohomology.

\hfill

\claim \label{_quotient_by_finite_homo_ratio_Claim_}
Let $f:\; X \arrow Y$ be a finite dominant morphism of complex varieties.
Assume that $X$ is smooth. Then $Y$ is weakly homology rational.
Moreover, if $Y$ is also normal, then $Y$ is homology rational.

{\bf Proof:} Well known. \endproof

\hfill

\definition
Let $\pi:\; X \arrow Y$ be a morphism of complex varieties.
Consider a stratification $\{C_i\}$ of $Y$, defined in such 
a way that the restriction of $f$ to the 
open strata $\pi^{-1}(C_i)$ is a locally trivial fibration. 
The map $f$ is called {\bf a semismall resolution
of $Y$} if $X$ is smooth, and for all $i$, the dimension
of the fibers of $\pi:\; \pi^{-1}(C_i) \arrow C_i$
is at most half the codimension of $C_i \subset Y$.

\hfill

\proposition \label{_cohomo_semismall_gene_Proposition_}
Let $\pi:\; X \arrow Y$ be a semismall resolution,
associated with the stratification $Y = \coprod C_i$.
Let $Y_i$ be the closed strata of the corresponding
stratification of $Y$, $Y_i = \bar{C_i}$,
and $X_i$ the corresponding closed strata of $X$,
$Y_i = \pi^{-1}(X_i)$. 
Consider the Goresky-MacPherson sheaf $V_i$ associated
with the sheaf $R^{l_i}\pi_* \C_{X_i}$
where $l_i=\frac{\codim_\C Y_i}{2}$
(for $\codim_\C Y_i$ odd, we put 
$V_i=0$) and $\C_{X_i}$ a constant 
sheaf on $X_i$. Assume that
$X$ is a K\"ahler manifold. 
Then $R^\bullet \pi_* \C_X$
is a direct sum of the perverse sheaves
$V_i$ shifted by $l_i$.

{\bf Proof:} In algebraic situation,
this is proven using the weight arguments
and $l$-adic cohomology (\cite{_Asterisque_100_}). 
To adapt this argument in K\"ahler situation,
one uses the mixed Hodge modules of M. Saito,
\cite{_Saito_}.
\endproof

\hfill

\theorem \label{_Hilbert_sche_semismal_cohomo_Theorem_}
\cite{_Gott_Sorg_}
Let $M$ be a complex surface,
$M^{[n]}$ its Hilbert scheme,
$M^{(n)}$ the symmetric power and
\[ \pi:\; M^{[n]}\arrow M^{(n)}\] the standard projection
map. Consider the stratification of $M^{(n)}$
by the diagonals $\Delta_{(\alpha)}$, parametrized
by the Young diagrams $\alpha$ (see 
\eqref{_Del-a_alpha_definition_Equation_}).
Then $\pi:\; M^{[n]}\arrow M^{(n)}$  is a semismall
resolution associated with this stratification.
Moreover, the sheaves $V_i$ of 
\ref{_cohomo_semismall_gene_Proposition_}
are constant sheaves $\C_{\Delta_{(\alpha)}}$.
\footnote{The varieties $\Delta_{(\alpha)}$ are normal
and obtained as quotients of smooth manifolds by group action.
Thus, all $\Delta_{(\alpha)}$ are
homology rational by \ref{_quotient_by_finite_homo_ratio_Claim_}.
Thus, the sheaves $\C_{\Delta_{(\alpha)}}$ are GM-type.}

{\bf Proof:} The map $\pi:\; M^{[n]}\arrow M^{(n)}$
is a semismall resulution, which is easy to check
by counting dimensions (see \ref{_holo_symple_semismall_Proposition_} 
for a conceptual proof). Now, the first assertion of 
\ref{_Hilbert_sche_semismal_cohomo_Theorem_}
is a straightforward application
of \ref{_cohomo_semismall_gene_Proposition_}. The second assertion
is much more subtle; see \cite{_Nakajima_} for details
and further reference. \endproof

\hfill

\corollary \label{_cohomo_of_Hilbert_explicit_Corollary_}
The $i$-th cohomology of $M^{[n]}$ are isomorphic to
\begin{equation} \label{_decompo_of_cohomo_Equation_}
  \bigoplus_\alpha 
  {\text{\Large \it H}}
  ^{\text{\large i}
  + \frac{\codim\Delta_{(\alpha)}}{2}}\left(\Delta_{(\alpha)}\right)
\end{equation}

{\bf Proof:} By \ref{_Hilbert_sche_semismal_cohomo_Theorem_},
\[
  R^\bullet\pi_* \C_{M^{[n]}} = \oplus \C_{\Delta_{(\alpha)}} 
  \left[ \frac{\codim\Delta_{(\alpha)}}{2}\right ],
\]
where $[\cdots]$ denotes the shift by that number.
\endproof

\subsection[Holomorphically symplectic manifolds and semismall
resolutions]{Holomorphically symplectic manifolds \\and semismall
resolutions}

\hfill

\definition
Let $\pi:\; X \arrow Y$ be a morphism of complex varieties,
and $\c Y$, $\c X$ be stratification of $Y$, $X$. Denote by
$Y_i$ the strata of $\c Y$.
We say that $\c Y$ and $\c X$ are compatible, if the preimages
$\pi^{-1}(Y_i)$ coinside with the strata $X_i$
of $\c X$, all the strata of $\c Y$
are non-singular, and the maps $\pi\restrict{X_i}:\; X_i \arrow Y_i$
are locally trivial fibrations. 

\hfill

\proposition \label{_holo_symple_semismall_Proposition_}
Let $\pi:\; X \arrow Y$ be a generically finite, dominant morphism
of complex varieties. Assume that $X$ is smooth and 
equipped with a holomorphically symplectic form
$\Theta_X$. Moreover, assume that there exists a Cartesian square
\[
\begin{CD}
\tilde X @>{\tilde \pi}>>  \tilde Y\\
@VV{\tilde \sigma}V @VV{\sigma}V \\
X @>{\pi}>>  Y
\end{CD}
\]
with finite dominant morphisms as 
vertical arrows and birational morphisms as horisontal arrows.
\footnote{By a {\bf birational morphism}
 we understand a morphism $\phi:\; X_1 \arrow X_2$
of complex varieties such that the inverse of $\phi$ is rational.}
Assume that $\tilde Y$ is a holomorphically symplectic manifold,
and the pullbacks of the holomorphic symplectic forms 
$\Theta_X$, $\Theta_{\tilde Y}$ via
${\tilde \pi}$ and $\tilde \sigma$ are equal on subvarieties,
in the sense of Subsection \ref{_symple_on_Hi_Subsection_}.
Assume, finally, that there exist compatible stratifications
$\{X_i\}$, $\{\tilde X_i\}$, $\{\tilde Y_i\}$ 
such that $\Theta_{\tilde Y}\restrict{\tilde Y_i}$ 
is non-degenerate in the generic
points of $\tilde Y_i$.
Then $\pi:\; X \arrow Y$ is a semismall resolution.

\hfill

{\bf Proof:}
Let $r(X_i)$ be the rank of the radical of $\Theta_X\restrict{X_i}$ in the
generic point of the stratum $X_i$. Similarly, let
$r(\tilde X_i)$ the rank of the radical of 
$\tilde \sigma^*\Theta_X\restrict{\tilde X_i}$ in the
generic point of the stratum $\tilde X_i$. Since
$\tilde \sigma$ is finite dominant, we have $r(\tilde X_i)= r(X_i)$.
By definition, $\tilde Y_i = \tilde\pi(\tilde X_i)$ is (generically) 
a non-degenerate symplectic subvariety of $\tilde Y$.
Since the forms $\Theta_{\tilde X}$ and 
$\tilde \pi^* \Theta_{\tilde Y}$ are equal on subvarieties,
and $\Theta_{\tilde Y}\restrict{\tilde Y_i}$ 
is non-degenerate in the generic
points of $\tilde Y_i$, we have 
\begin{equation}\label{_r_X_i_=_dim_fib_Equation_}
   r(\tilde X_i) = \dim_\C \left(\tilde\pi^{-1}(y)\right),
\end{equation}
for $y\in \tilde Y_i$ a generic point.
Let $w(X_i)$ be the number $\dim(X_i) - r(X_i)$.
The following linear-algebraic claim immediately implies that
\begin{equation}\label{_codim_Y_i_>=_2r_X_i_Equation_}
   \dim_\C X - w(X_i) \geq 2 r(X_i)
\end{equation}
\endproof

\hfill

\claim 
Let $W$ be a symplectic vector space, $\Theta$ the symplectic
form, $V\subset W$ a subspace, $r(V)$ the rank of the radical 
$\Theta\restrict V$ and $w(V):= \dim V - r(V)$.
Then $\dim W - w(V) \geq 2 r(V)$.

{\bf Proof:} Clear. \endproof

\hfill

Comparing \eqref{_r_X_i_=_dim_fib_Equation_}
and \eqref{_codim_Y_i_>=_2r_X_i_Equation_},
we obtain that 
\[ \codim_\C Y_i \geq 2 \dim_\C \left(\tilde\pi^{-1}(y)\right), \]
for $y\in \tilde Y_i$ a generic point. Finally, since 
$\sigma:\; \tilde Y \arrow Y$ is finite dominant, we have
$\dim_\C \left(\tilde\pi^{-1}(y)\right) = 
\dim_\C \left(\pi^{-1}(\sigma(y))\right)$.
Thus, $\codim_\C Y_i \geq 2 \dim_\C \left(\pi^{-1}(y)\right)$
for $y\in Y_i$ a generic point. This finishes the proof
of \ref{_holo_symple_semismall_Proposition_}. \endproof

\hfill

\corollary \label{_triana_speci_semismall_Corollary_}
Let $M$ be a complex K3 surface or a compact complex 2-dimensional torus, 
$M^{[n]}$ its Hilbert scheme and $M^{(n)}$ the symmetric power
of $M$, with $\pi:\; M^{[n]}\arrow M^{(n)}$
being the standard map. Consider an arbitrary hyperk\"ahler
structure $\c H$ on $M^{[n]}$ compatible with the complex
structure. Let $Z\subset M^{[n]}$ be a subvariety
which is trianalytic with respect to
$\c H$, and $n:\; X \arrow Z$
be the desingularization of $Z$. Assume that
$M$ is Mumford-Tate generic
with respect to some hyperk\"ahler structure.
Then $\pi \circ n :\; X \arrow Y$ is
a semismall resolution of $Y:=\pi(Z)$. 

\hfill

{\bf Proof:} Assume that $Z$ is irreducible. 
Since the desingularization
$X$ is hyperk\"ahler, this manifold is holomorphically
symplectic, and the holomorphic symplectic form $\Theta_{X}$
on $X$ is obtained as a pullback of the holomorphic
symplectic form $\Theta_{M^{[n]}}$ on $M^{[n]}$.
To simplify notations, we denote $\pi \circ n$ by $\pi$.
Let 
\begin{equation}\label{_commu_squa_with_subva_Equation_}
\begin{CD}
\tilde X @>{\tilde \pi}>>  \tilde Y\\
@VV{\tilde \sigma}V @VV{\sigma}V \\
X @>{\pi}>>  Y
\end{CD}
\end{equation}
be the Cartesian square,
with $\tilde Y$ obtained as an irreducible component
of the preimage $\sigma^{-1}(Y)\subset M^n$. We 
intend to show that the
square \eqref{_commu_squa_with_subva_Equation_}
satisfies the assumptions of 
\ref{_holo_symple_semismall_Proposition_}.
For each morphism of complex varieties, there
exists a stratification, compatible with this
morphism. Take a set of compatible stratifications $\{X_i\}$,
$\{\tilde Y_i\}$, $\{\tilde X_i\}$.
By \ref{_hyperkae_embeddings_Corollary_}, 
any stratification of $\tilde Y$ consists of 
trianalytic subvarieties because all closed complex subvarieties
of $M^n$ are trianalytic. 
Applying \ref{_holo_symple_semismall_Proposition_}
to the map $\pi:\; X \arrow Y$ and the Cartesian square
\eqref{_commu_squa_with_subva_Equation_},
we immediately obtain \ref{_triana_speci_semismall_Corollary_}.
\endproof


\section{Universal subvarieties of the Hilbert scheme}
\label{_unive_subva_Section_}


The Sections \ref{_unive_subva_Section_}--\ref{_triana_unive_subva_Section_} 
are independent from the rest of this paper. 
The only result of 
Sections \ref{_unive_subva_Section_}--\ref{_triana_unive_subva_Section_} 
that we use is \ref{_one-to-one_triana_Corollary_}.

Let $M$ be a smooth complex surface, $M^{[n]}$ its
Hilbert scheme. An automorphism $\gamma$ of $M$ gives an automorphism
$\gamma^{[n]}$ of $M^{[n]}$; similarly, an infinitesimal automorphism
of $M$ (that is, a holomorphic vector field) gives
an infinitesimal automorphism of $M^{[n]}$. 

\hfill

\definition \label{_Unive_subva_Definition_}
Let $M$ be a surface, $M^{[n]}$ its
Hilbert scheme and $X\subset M^{[n]}$ a closed complex subvariety.
Then $X$ is called {\bf universal} if for all
open $U\subset M$, and all global or infinitesimal
automorphisms $\gamma\in  \Gamma(T(M))$,
the subvariety $X_U := X \cap U^{[n]}$ is preserved by
$\gamma^{[n]}$.

\hfill

The universal subvarieties are described more explicitly
in the following subsection

\subsection{Young diagrams and universal subvarieties
of the Hilbert scheme}
\label{_Young_Hilb_Subsection_}

Let $M$ be a smooth surface, $M^{(n)}$ its
symmetric power, $M^{[n]}$ its Hilbert scheme
and $\pi:\; M^{[n]}\arrow M^{(n)}$ the natural map.
For a Young diagram 
\[ \alpha = (n_1\geq n_2\geq ... \geq n_k), \ \ \sum n_i =n, \]
we defined a subvariety 
$\Delta_{(\alpha)}\subset M^{(n)}$
\eqref{_Del-a_alpha_definition_Equation_}. Let 
$\Delta_{[\alpha]}:= \pi^{-1}(\Delta_{(\alpha)})$
the the corresponding subvariety in $M^{[n]}$.

Let $a$ be the general point of $\Delta_{(\alpha)}$, 
i. e. the one satisfying
\begin{equation}\label{_gene_of_Delta_alpha_Equation_}
\begin{split}
a = &\sigma(a_1, ..., a_n),\ \ \text{where} \ \  
         a_1 = a_2 = ... = a_{n_1}  \\&
               a_{n_1+1} = a_{n_1+2} = ... = a_{n_1+n_2} \  
               ..., \\ &
               a_{\sum_{i=1}^{k-1} n_i+1} = 
               a_{\sum_{i=1}^{k-1} n_i+2} = ... = a_n \\
& 
\text{and the points $a_1$, $a_{n_1+1}$, ..., 
$a_{\left(\sum_{i=1}^{k-1} n_i\right)+1}$ are pairwise unequal}
\end{split}
\end{equation}
Let $F_\alpha(a):= \pi^{-1}(a) \subset \Delta_{[\alpha]}$
be the general fiber of the projection
$\pi:\; \Delta_{[\alpha]} \arrow \Delta_{(\alpha)}$.
By definition, $F_\alpha(a)$ parametrizes $0$-dimensional subschemes
$S\subset M$ with $Sup(S)= \{a_i\}$ and 
prescribed multiplicites \[ \operatorname{length}_{a_i}S=n_i. \]
Let $\hat \calo_{a_i}$ be the adic completion of 
$\calo_M$ at $a_i$, and $G_{a_i} := \Aut(\hat \calo_{a_i})$.
Clearly, the group $G_a:= \prod_i G_{a_i}$ acts naturally
on $F_\alpha(a)$. We are interested in $G_a$-invariant subvarieties
of $F_\alpha(a)$.

\hfill

\lemma \label{_inva_subva_F_alpha_identifi_Lemma_}
Let $\alpha$ be a Young diagram, $\Delta_{(\alpha)}$ the corresponding
subvariety of $M^{(n)}$ and $a, b$ the points of $\Delta_{(\alpha)}$
satisfying \eqref{_gene_of_Delta_alpha_Equation_}. Let
$F_\alpha(a)$, $F_\alpha(b)$ be the corresponding fibers of
$\pi:\; \Delta_{[\alpha]} \arrow \Delta_{(\alpha)}$. Consider the
groups $G_a$, $G_b$ acting on $F_\alpha(a)$, $F_\alpha(b)$ as
above. Then 
\begin{description}
\item[(i)] there exist a canonical bijective
correspondence $\theta$ between $G_a$-invariant subvarieties
in $F_\alpha(a)$ and $G_b$-invariant subvarieties
in $F_\alpha(b)$.
\item[(ii)]
For any complex automorphism $\gamma\; M \arrow M$
such that $\gamma(a)=b$, the corresponding map
$\gamma:\; F_\alpha(a)\arrow F_\alpha(b)$ maps
$G_a$-invariant subvarieties
of $F_\alpha(a)$ to $G_b$-invariant subvarieties
of $F_\alpha(b)$ and induces 
the correspondence $\theta$.
\end{description}

{\bf Proof:} Let $(a_1, ... a_n)$, $(b_1, ... b_n)$
be the points of $M^n$ satisfying
\eqref{_gene_of_Delta_alpha_Equation_}, 
such that $a = \sigma(b_1, ... b_n)$,
 $b=\sigma(b_1, ... b_n)$. Let $U$ be an open subset of $M$ 
containing $a_i$, $b_i$, $i = 1, ... , n $. Let 
$\gamma:\; U \arrow U$ be a complex automorphism of $U$ such that
$\gamma(a_i) = b_i$. Since $a$, $b$ satisfy 
\eqref{_gene_of_Delta_alpha_Equation_}, for 
an appropriate choice of $U$, such $\gamma$ 
always exists. Clearly, $\gamma$ identifies $F_\alpha(a)$
and $F_\alpha(b)$. This identification is {\it not} unique,
since it depends on the choice of $\gamma$, but
every two such identifications differ by a twist by
$G_a$, $G_b$. This proves
\ref{_inva_subva_F_alpha_identifi_Lemma_}.
\endproof

\hfill

By \ref{_inva_subva_F_alpha_identifi_Lemma_}, the set of
$G_a$-invariant subvarieties of $F_\alpha(a)$ is independent from
$a$. Denote this set by $\c Z_\alpha$. For each 
$Y\in \c Z_\alpha$, and a generic point 
$a\in \Delta_{(\alpha)}$, denote the corresponding
subvariety of $F_\alpha(a)$ by $Y(a)$. Let $\c Z_\alpha(Y)$
be the union of $Y(a)$ for all $a\in \Delta_{(\alpha)}$
satisfying \eqref{_gene_of_Delta_alpha_Equation_}.

\hfill

\theorem\label{_inva_subva_from_Young_Theorem_}
Let $\alpha$ be a Young diagram, 
$\Delta_{(\alpha)}\subset M^{(n)}$ the corresponding
diagonal in $M^{(n)}$ and $a\in \Delta_{(\alpha)}$
a general point (that is, one 
satisfying \eqref{_gene_of_Delta_alpha_Equation_}).
Let $Y\in \c Z_\alpha$ be a $G_a$-invariant subvariety of 
$F_\alpha(a)= \pi^{-1}(a)\subset M^{[\alpha]}$,
and $\c Z_\alpha(Y)$ the corresponding subvariety of 
$M^{[\alpha]}$. Then $\c Z_\alpha(Y)$
is a universal subvariety of
$M^{[\alpha]}$, in the sense of \ref{_Unive_subva_Definition_}.
Moreover, all irreducible universal subvarieties of
$M^{[\alpha]}$ can be obtained  this way.

\hfill

{\bf Proof:} The statement of 
\ref{_inva_subva_from_Young_Theorem_} is local by
$M$. Thus, to prove that $\c Z_\alpha(Y)$ is preserved
by infinitesimal automorphisms, it suffices to show that
$\c Z_\alpha(Y)$ is preserved by all global automorphisms 
of $M$. Let $\gamma:\; M \arrow M$ be an automorphism.
Denote by $\Delta_{(\alpha)}^\circ\subset \Delta_{(\alpha)}$
the set of all $a$ satisfying \eqref{_gene_of_Delta_alpha_Equation_}.
Clearly, $\gamma$ preserves 
\[ \Delta_{[\alpha]}^\circ:= 
   \pi^{-1}\left(\Delta_{(\alpha)}^\circ\right).
\] 
To show that $\gamma$ preserves $\c Z_\alpha(Y)$,
it suffices to prove that, for all $a, b \in 
\Delta_{(\alpha)}^\circ$, $\gamma(a)=b$, the automorphism
$\gamma$ maps $F_\alpha(a)$ to $F_\alpha(b)$. This is
\ref{_inva_subva_F_alpha_identifi_Lemma_} (ii). 
We obtained that $\c Z_\alpha(Y)$ is universal.

\hfill

Let $X$ be an irreducible universal subvariety in
$M^{[n]}$. Then $\pi(X)\subset M^{(n)}$ is preserved by
the automorphisms of $M^{(n)}$ coming 
from $M$. For $M$ Stein, the 
only subvarieties of $M^{(n)}$ preserved by infinitesimal
automorphisms are unions of diagonals.
Since $X$ is irreducible, so is $\pi(X)$. We obtain that 
$\pi(X)$ is a diagonal $\Delta_{(\alpha)}$ corresponding to 
some Young diagram $\alpha$. It remains to prove that
$X\cap F_\alpha(a)$ is $G_a$-invariant, for all
$a\in \Delta_{(\alpha)}^\circ$. This is clear, because
infinitesimal automorphisms  of $M$ fixing $\{a_i\}$ generate the
group $G_a =\prod_i \Aut(\hat \calo_{a_i})$, and since
$X$ is invariant under such automorphisms, 
$X\cap F_\alpha(a)$ is $G_a$-invariant. 
\ref{_inva_subva_from_Young_Theorem_} is proven.
\endproof

\subsection{Universal subvarieties of relative dimension 0}

\definition
Let $M$ be a smooth complex surface, $M^{[n]}$ its Hilbert
scheme, $\alpha$ a Young diagram corresponding to
a diagonal $\Delta_{(\alpha)} \subset M^{(n)}$. Let
$a\in \Delta_{(\alpha)}$ be a general point,
and $F_\alpha(a):= \pi^{-1}(a)$ the corresponding fiber
of $\pi:\; M^{[n]} \arrow M^{(n)}$.
Consider a $G_\alpha(a)$-invariant subvariety 
$Y$ of $F_\alpha(a)$. Let $Z\subset M^{[n]}$ be a corresponding 
universal subvariety, $Z = Z_\alpha(Y)$ 
(\ref{_inva_subva_from_Young_Theorem_}).
Then the {\bf relative dimension} of $Z$ is the dimension of $Y$.

\hfill

In this subsection, we classify the universal subvarieties
of relative dimension 0. 

\hfill

Let $\alpha = (n_1\geq n_2\geq ... \geq n_k)$ 
be a Young diagram, $\sum n_i =n$. 
Clearly, 
\begin{equation}\label{_F_alpha_via_F_0_Equation_}
F_\alpha(a) \cong F_0(n_1) \times F_0(n_2) \times ... ,
\end{equation}
where $F_0(i)$ is the classifying space of $0$-dimensional
subschemes of length $i$ in $\C^2$ with support in $0\in \C^2$.
Let $G_0 = \Aut(\C[[x,y]])$ be the group of automorphisms of the
ring of formal series, acting on $F_0(i)$. By 
\eqref{_F_alpha_via_F_0_Equation_},
the $k$-th power of $G_0$ acts on $F_\alpha(a)$.
This gives an isomorphism $G_0^k \cong G_\alpha(a)$.

Assume that
$n_i = \frac{m_i\cdot (m_i+1)}{2}$, for some 
positive integer $m_i$. Consider a $G_0$-invariant point
$s_i\in F_0(n_i)$, given by
\begin{equation} \label{_inva_point_quotie_by_power_Equation_}
   s_{m_i} = \C[[x,y]] / {\goth m}^{m_i},
\end{equation}
where $\goth m\subset \C[[x,y]]$ is the maximal 
ideal generated by $x$ and $y$.
Let $\{s_1\} \times\{s_2\} \times ... \times\{s_k\}$
be the $G_0^k$-invariant point of $\prod F_0(n_i)$.
Using the isomorphism \eqref{_F_alpha_via_F_0_Equation_},
we obtain a $G_\alpha(a)$-invariant point $a$ of $F_\alpha(a)$.
Denote by $\c X_\alpha$ the corresponding universal subvariety
of $M^{[n]}$. It has relative dimension 0. The aim
of this subsection is to show that
all universal subvarieties of relative dimension $0$ are
obtained this way.

\hfill

\proposition \label{_unive_subva_rela_dime_0_Proposition_}
Let $X\subset M^{[n]}$ be a universal subvariety of relative dimension
$0$. Then where exists a Young diagram
\[ \alpha = (n_1\geq n_2\geq ... \geq n_k), \ \ \sum n_i =n, \]
and positive integers
$m_1, ..., m_k$, such that $n_i = \frac{m_i\cdot (m_i+1)}{2}$,
and $X = \c X_\alpha$.

\hfill

{\bf Proof:} Let $a$ be a general point of
$\Delta_\alpha$, and $s\in F_\alpha(a)$ a point
of the zero-dimensional variety $F_\alpha(a)\cap X$.
Consider the varieties $F_0(i)$ defined
above, and the action of $G_0 = \Aut(\C[[x,y]])$
on $F_0(i)$. Let $x_i\in F_0(n_i)$ be the points
of $F_0(k)$, such that under an isomorphism
\eqref{_F_alpha_via_F_0_Equation_}, $s$ corresponds
to $\{x_1\} \times \{x_2\} \times ... \times \{x_k\}$.
Then the points $x_i$ are $G_0$-invariant.
To finish the proof of \ref{_unive_subva_rela_dime_0_Proposition_},
it remains to prove the following lemma.

\hfill

\lemma \label{_G_0_inva_poins_in_F_0_Lemma_}
Let $s\in F_0(i)$ be a $G_0$-invariant point. Then 
$i= \frac{j\cdot(j+1)}{2}$ and $s$ is given by
\eqref{_inva_point_quotie_by_power_Equation_}.

{\bf Proof:} The group $GL(2, \C)$ acts on $\C[[x,y]]$ 
by automorphisms. Clearly, this $GL(2, \C)$-action is factorized
through the natural action of \[ G_0 = \Aut(\C[[x,y]]). \]
We show (\ref{_GL(2)_inva_ideals_in_series_Sublemma_} below) that
all $GL(2, \C)$-invariant ideals in $\C[[x,y]]$ are powers of the
maximal ideal. Since $x = \C[[x,y]]/ I$ for some
$G_0$-, and hence, $GL(2, \C)$-invariant ideal of
$\C[[x,y]]$, this will finish the proof of
\ref{_G_0_inva_poins_in_F_0_Lemma_}.
We reduced \ref{_unive_subva_rela_dime_0_Proposition_}
to the following result.

\hfill

\sublemma \label{_GL(2)_inva_ideals_in_series_Sublemma_}
Consider the natural action of $GL(2, \C)$ on \[ A= \C[[x,y]].\]
Let $I$ be a proper $GL(2, \C)$-invariant
ideal in $A$. Then $I$ is a power of the
maximal ideal. 

{\bf Proof:} 
Consider the $GL(2)$-invariant filtration
\[ 
  A_0 \subset A_0 \oplus A_1 \subset A_0\oplus A_1 \oplus
A_2 \subset ... \subset A
\]
where $A_i\subset A$ consists of homogeneous polynomials
of degree $i$.
Let $l$ be the minimal number for which $I\cap A_l \neq 0$. 
Since $I$ and $A_l$
are $GL(2)$-invariant, the intersection $I \cap A_l$
is also $GL(2)$-invariant. 
The space $A_l$ is an irreducible representation
of $GL(2)$, and thus, $I \supset A_l$.
Therefore, $I = A_l \cdot A$, and $I$ is $l$-th
power of the maximal ideal. This finishes the
proof of \ref{_GL(2)_inva_ideals_in_series_Sublemma_},
\ref{_G_0_inva_poins_in_F_0_Lemma_}, and
\ref{_unive_subva_rela_dime_0_Proposition_}.
\endproof



\section{Subvarieties of $M^{[n]}$ which are generically 
finite over $M^{(n)}$, for $M$ a generic K3 surface}
\label{_subva_gene_fini_Section_}


Throughout this section, $M$ is a smooth complex surface,
$M^{[n]}$ the Hilbert scheme of $M$,
$M^{(n)}$ the $n$-th symmetric power of $M$ and
$\pi:\; M^{[n]}\arrow M^{(n)}$ the natural map.

Let $f:\; X \arrow Y$ be a morphism of complex varieties.
We say that $f$ is {\bf generically finite} if
there exist an open dense subset $X_0\subset X$ 
such that the map $f\restrict{X_0}:\; X_0 \arrow f(X_0)$
is finite. The morphism $f$ is called {\bf generically
one-to-one} if there exist an open dense subset $X_0\subset X$ 
such that the map $f\restrict{X_0}:\; X_0 \arrow f(X_0)$
is an isomorphism.

\hfill

The main result of this section is the following theorem.

\hfill

\theorem\label{_gene_fini_universa_Theorem_}
Let $M$ be a complex K3 surface.
Assume that $M$ admits a hyperk\"ahler
structure $\c H$ such that $M$ is Mumford-Tate
generic with respect to $\c H$ (\ref{_generic_manifolds_Definition_}).
Let $X\subset M^{[n]}$ be an irreducible complex 
analytic subvariety such that the restriction of
$\pi$ to $X$ is generically finite. Assume that
there exists a Young diagram $\alpha$ such that
the subvariety $\pi(X)\subset M^{(n)}$ coinsides 
with $\Delta_{(\alpha)}$.
Then $X$ is a universal subvariety (\ref{_Unive_subva_Definition_})
 of $M^{[n]}$.

\hfill

\remark
The relative dimension of the universal subvariety 
$X\subset M^{[n]}$ is zero, because $\pi\restrict X$ is generically
finite. Thus, \ref{_unive_subva_rela_dime_0_Proposition_}
can be applied to this situation. We obtain that,
under assumptions of \ref{_gene_fini_universa_Theorem_},
$\pi\restrict X :\; X \arrow \Delta_{(\alpha)}$
is generically one-to-one. 

\hfill

The proof of \ref{_gene_fini_universa_Theorem_} takes
the rest of this section.

\subsection{Fibrations arising from the Hilbert scheme}

We work in assumptions of \ref{_gene_fini_universa_Theorem_}.
Let $\Delta^\circ_{(\alpha)}\subset \Delta_{(\alpha)}$ be the 
set of general points of $\Delta_{(\alpha)}$, defined by
\eqref{_gene_of_Delta_alpha_Equation_}. Consider the
fibration $\pi:\; \Delta^\circ_{[\alpha]}\arrow 
\Delta^\circ_{(\alpha)}$, where 
$\Delta^\circ_{[\alpha]}= \pi^{-1}(\Delta^\circ_{(\alpha)})$.

Let $M^{(l)}_\circ$ be the $M^{(l)}$ with all diagonals deleted:
\[ M^{(l)}_\circ = 
   \bigg\{ (x_1 , ... x_l) \in M^{(l)}\ \ |\ \ x_i\neq x_j 
    \text{\ \ for all \ \ } i\neq j
   \bigg\}
\]
We write $\alpha = (n_1 \geq n_2 \geq ... \geq n_k )$
as 
\begin{equation}\label{_n'_i_defi_Equation_}
\begin{split}
  \alpha = &\left(n_1= n_2 = ... =n_{n_1'} > 
   n_{n'_1+1} = ... = n_{n'_1 +n'_2} > 
\vphantom{ n_{\sum_{i=1}^{k'-1} n_i' -1} } \right.
   ... \\& \left. ... > n_{\sum_{i=1}^{k'-1} n_i' +1} = ... = 
    n_{\sum_{i=1}^{k'-1} n_i' -1} = n_{\sum_{i=1}^{k'} n_i'}\right),
\end{split}
\end{equation}
where $\sum_{i=1}^{k'} n_i' =k$. 

\hfill

\claim\label{_Delta_without_diags_product_Claim_}
The manifold $\Delta^\circ_{(\alpha)}$ is naturally isomorphic
to $\prod_i M^{(n_i')}_\circ$. 

{\bf Proof:} Clear. \endproof

\hfill

Let $D^\circ_{(\alpha)}$
be the universal covering of $\Delta^\circ_{(\alpha)}$.
{}From \ref{_Delta_without_diags_product_Claim_} it is clear
that 
\begin{equation}\label{_D_alpha_is_M^K'_Equation_}
   D^\circ_{(\alpha)} = \prod_i M^{n_i'}_\circ\subset M^{k'},
\end{equation}
where $M^{n_i'}_\circ$ is $M^{n_i'}$ without diagonals.
We define $D^\circ_{[\alpha]}$ as a fibered product,
in such a way that the square
\begin{equation}\label{_D_alpha_Cartesian_defi_Equation_}
\begin{CD}
D^\circ_{[\alpha]} @>>> \Delta^\circ_{[\alpha]}\\
@VVpV @VV\pi V \\
D^\circ_{(\alpha)} @>>> \Delta^\circ_{(\alpha)}
\end{CD}
\end{equation}
is Cartesian. The 
map $D^\circ_{[\alpha]}\stackrel p \arrow D^\circ_{(\alpha)}$
is a locally trivial fibration. We determine the fibers
of $p$ in terms of the isomorphism 
\eqref{_D_alpha_is_M^K'_Equation_} as follows.

\hfill

Consider the vector bunlde $J^i(M)$ over $M$, with the fibers
$J^i(M)\restrict x = \calo_M/ {\goth m}_x^i$, where
${\goth m}_x$ is the maximal ideal on $\calo_M$ corresponding to $x$.
Clearly, the bundle $J^i(M)$ has a natural ring structure.
Let $G^i(M)$ be the fibration over $M$ with the fibers
$G^i(M)\restrict x$ classifying the ideals
$I\subset J^i(M)$ of codimension $i$.
Consider again the equation \eqref{_n'_i_defi_Equation_}. Let 
$N:\; \{ 1 , .. k'\} \arrow \Z^+$ be the map
\[ l \arrow n_{\sum_{i=1}^{l-1}}, 
\]
i. e., $1$ is mapped to $n_1$, $2$ to the biggest value of $n_i$ not
equal to $n_1$, $3$ to the third biggest, etc.

For a locally trivial fibrations $Y_1$, $Y_2$
over $X_1$, $X_2$, we denote the external product
by $Y_1 \newboxtimes Y_2$. This is a fibration over 
$X_1\times X_2$, with the fibers which are products
of fibers of $Y_1$, $Y_2$. The iterations of $\newboxtimes$
(for three or more fibrations) are defined in the same 
spirit.

\hfill

\claim
Under the isomorphism \eqref{_D_alpha_is_M^K'_Equation_},
the locally trivial fibration
$p:\; D^\circ_{[\alpha]}\arrow D^\circ_{(\alpha)}$
is isomorphic to the fibration 
\[
   \newboxtimes\limits_{i=1}^{k'} G^{N(i)}(M) \restrict {D^\circ_{(\alpha)}}
\]
over $D^\circ_{(\alpha)} \subset M^{k'}$.

{\bf Proof:} Clear. \endproof

\hfill

Let $D_{[\alpha]}\arrow D_{(\alpha)}$ be the fibration
$\newboxtimes\limits_{i=1}^{k'} G^{N(i)}(M) \arrow M^{k'}$.
We consider $D^\circ_{[\alpha]}$, $D^\circ_{(\alpha)}$
as open subsets in $D_{[\alpha]}$, $D_{(\alpha)}$.

\hfill

Let $X\subset M^{[n]}$ be a closed subvariety,
$\pi(X) = \Delta_{(\alpha)}$, and
$\pi:\; X \arrow \Delta_{(\alpha)}$ generically finite.
Consider $X\cap \Delta^\circ_{[\alpha]}$ as a closed 
subvariety of $\Delta^\circ_{[\alpha]}$. Let 
$\tilde X$ be an irreducible component of
$n^{-1} (X) \subset D^\circ_{[\alpha]}$,
where $n:\; D^\circ_{[\alpha]}\arrow \Delta^\circ_{[\alpha]}$
is the horisontal arrow of 
\eqref{_D_alpha_Cartesian_defi_Equation_}.
Clearly, the closure of $\tilde X$ in $D_{[\alpha]}$
is a closed complex subvariety of $D_{[\alpha]}$.
We denote this subvariety by $Z$. By construction,
$Z$ is irreducible (it is an image of an irreducible
variety) and generically finite
over $D_{(\alpha)} = M^{k'}$.

\hfill

Consider the fibration $G^m(M) \arrow M$ constructed above.
Assume that $m = \frac{l\cdot (l+1)}{2}$ for a positive 
integer $l$. Then the fibration $G^m(M) \arrow M$ has
a canonical section $s:\; M \arrow G^m(M)$, defined 
by $x \arrow \calo_M /{\goth m}_x^l$, where
${\goth m}_x\subset \calo_M$ is the maximal ideal 
corresponding to $x$.

\hfill

\proposition\label{_subvarie_of_G^i_Proposition_}
Let $M$ be a complex K3 surface.
Assume that $M$ admits a hyperk\"ahler structure $\c H$
such that $M$ is generic with respect to $\c H$.
Let $Y\subset G^m(M)$ be
a closed irreducible subvariety of the total space of the
fibration  $G^m(M) \stackrel p \arrow M$. 
Assume that $Y$ is generically finite over $M$.
Then $m = \frac{l\cdot (l+1)}{2}$ for some positive 
integer $l$, and $Y$ is the image of the natural map
$s:\; M \arrow G^m(M)$ constructed above.

\hfill

\ref{_subvarie_of_G^i_Proposition_} is proven
in Subsection \ref{_fibra_ove_K3_Subsection_}. 
Presently, we are going to 
explain how \ref{_subvarie_of_G^i_Proposition_} implies 
\ref{_gene_fini_universa_Theorem_}.

\hfill

Consider the map $p:\; D_{[\alpha]} \arrow M^{k'}$,
and the closed subvariety $Z\subset D_{[\alpha]}$
constructed from $X$ as above. 
Let $(m_1, ... , m_{k'-1}) \in M^{k'-1}$ be a point such that
the map
\[ p:\; Z \cap p^{-1}(\{(m_1, ... , m_{k'-1})\} \times M)
\arrow \{(m_1, ... , m_{k'-1})\} \times M)
\]
is generically finite.
The set $S$ of such $(m_1, ... , m_{k'-1})$ is open and dense in
$M^{k'-1}$. Let $\Psi_i:\;  D_{[\alpha]} \arrow G^{N(i)}(M)$
be the natural projection to the $i$-th component of the
product $D_{[\alpha]} = \newboxtimes\limits_{i=1}^{k'} G^{N(i)}(M)$.
By \ref{_subvarie_of_G^i_Proposition_},
$N(k')= \frac{l\cdot(l-1)}{2}$ and the subvariety
\[ \Psi(Z \cap p^{-1}(\{(m_1, ... , m_{k'-1})\} \times M))\]
coinsides with image of the map \[ s_{k'}:\; M \arrow G^{N(k')}(M).\]
Since $Z$ is irreducible,
\[ \Psi_{k'}(Z \cap p^{-1}(\{(m_1, ... , m_{k'-1})\} \times M))
   \subset G^{N(k')}(M)
\]
is independent from the choice of $(m_1, ... , m_{k'-1})\in S$.
Therefore, $\Psi_{k'}(Z) = \im(s_{k'})$. A similar argument shows
that $\Psi_{i}(Z) = \im(s_{i})$, for all $i= 1, ... , k'$.
Thus, $Z$ is an image of the section of the map
$p:\; D_{[\alpha]} \arrow M^{k'}$ given by
$\newboxtimes\limits_{i=1}^{k'} s_i$. This implies
\ref{_gene_fini_universa_Theorem_}. We reduced
\ref{_gene_fini_universa_Theorem_} to 
\ref{_subvarie_of_G^i_Proposition_}.

\subsection{Projectivization of stable bundles}

Let $M$ be a compact K\"ahler manifold. We
understand stability of holomorphic vector
bundles over $M$ in the sense of Mumford--Takemoto 
(\ref{_degree,slope_destabilising_Definition_}).
A polystable bundle is a direct sum of stable bundles
of the same slope. Let $V$ be a polystable 
 bundle, and ${\Bbb P} V$ its projectivization.
Consider the unique Yang-Mills connection on $V$
(\ref{Yang-Mills_Definition_}). 
This gives a natural connection $\nabla_V$
on the fibration ${\Bbb P} V\arrow M$.

\hfill

\proposition\label{_subvarie_of_PV_Proposition_}
Let $M$ be a compact complex simply connected 
manifold of hyperk\"ahler type.
Assume that $M$ admits a hyperk\"ahler structure $\c H$
such that $M$ is generic with respect to $\c H$.
Consider $M$ as a K\"ahler manifold, with 
the K\"ahler metric induced from $\c H$.
Let $V$ be a polystable
bunlde over $M$, and ${\Bbb P} V\stackrel \pi\arrow M$ 
its projectivization. Consider a closed irreducible 
subvariety  $X\subset {\Bbb P} V$ such that
$\pi(X) = M$. Then $X$ is preserved by the
connection $\nabla_V$ in ${\Bbb P} V$.

\hfill

{\bf Proof:} Let $x\in M$ be a point of $M$ such 
that in a neighbourhood $U\subset M$ 
of $x$, the projection $\pi:\; X \arrow M$
is a locally trivial fibration. Assume that
$U$ is open and dense in $M$. Let $X_x$ be the fiber of
$\pi:\; X \arrow M$ in $x$, and $V_x:= V\restrict x$ 
the fiber of $V$. Consider the Hilbert scheme $H$ classifying
the subvarieties $Y\subset {\Bbb P} V_x$
with the same Hilbert polynomial as $Y$.
Then $H$ can be naturally embedded to the projectivization
of a linear space $W_x$, where $W_x$ is a certain tensor power of
$V_x$, depending on the Hilbert polynomial of $X_x$.
Consider the corresponding bundle $W$, which is 
related to $V$ in the same way as $W_x$ to $V_x$.
Then, $W$ is a tensor power of $V$, and hence, $W$ is
equipped with a unique Yang-Mills connection. 
Consider the corresponding connection $\nabla_W$ on
the projectivization ${\Bbb P}W$. 
Let $X_0$ denote $\pi^{-1}(U) \cap X$.
The locally trivial fibration
$\pi\restrict{X_0}:\; X_0\arrow U$ gives a section 
$s$ of ${\Bbb P} W\restrict{X_0}$. To prove that
$X$ is preserved by the
connection $\nabla_V$ in ${\Bbb P} V$, it suffices to show that
$X_0$ is preserved by $\nabla_V$, or that
$\im s$ is preserved by $\nabla_W$.
This is implied by the following lemma,
which finishes the proof of 
\ref{_subvarie_of_PV_Proposition_}.

\hfill

\lemma\label{_secti_preserved_by_nabla_Lemma_}
In assumptions of \ref{_subvarie_of_PV_Proposition_}.
let $U\subset M$ be a dense open set, such that
$\pi\restrict{X_0}:\; X_0\arrow U$ is an isomorphism,
where \[ X_0 = \pi^{-1}(U) \cap X\subset {\Bbb P}V. \]
Then $\nabla_V$ preserves $X$.

\hfill

{\bf Proof:} Since $M$ is generic with respect to $\c H$,
all its complex subvarieties have complex codimension at least 
2. Thus, we may assume that the complement $M \backslash U$
is a complex subvariety of codimension at least 2.

Consider the restriction $V\restrict U$.
Then $X_0$ gives a one-dimensional subbundle $L$ of $V\restrict U$.
Let $V'= i_* L \subset i_* V\restrict U$ 
be the direct image of $L$ under the
embedding $U\stackrel i \hookrightarrow M$. Since $M \backslash U$
is a complex subvariety of codimension at least 2, 
the natural map $V \arrow i_* V\restrict U$ 
is an isomorphism. Therefore, $V'$ is a coherent subsheaf in
$V$. To prove \ref{_secti_preserved_by_nabla_Lemma_}
it suffices to show that $V'$ is preserved by the connection.
Since $M$ is generic with respect to $\c H$, all integer 
$(1,1)$-classes of cohomology have degree 0 
(\ref{_Lambda_of_inva_forms_zero_Lemma_}).
Therefore, $\slope(V')=\slope(V)=0$ and $V'$ is a destabilising
subsheaf of $V$. Since $V$ is polystable, this implies that
$V'$ is a direct summand of $V$, and the Yang-Mills connection
in $V$ preserves the decomposition $V = V' \oplus {V'}^{\bot}$, where
${V'}^{\bot}$ is the orthogonal complement of $V'$ with respect to
any Yang-Mills metric on $V$. \ref{_secti_preserved_by_nabla_Lemma_}
is proven. This finishes the proof of 
\ref{_subvarie_of_PV_Proposition_}. \endproof

\hfill

\corollary
In assumptions of \ref{_subvarie_of_PV_Proposition_},
let $X$ be generically finite over $M$. Assume that
$M$ is simply connected. Then $\pi:\; X\arrow M$ is an isomorphism.

{\bf Proof:} Since $X$ is preserved by the connection, the map
$\pi:\; X\arrow M$ is a finite covering. Since $M$ is simply 
connected, and $X$ is irreducible, $\pi:\; X\arrow M$ is one-to-one.
\endproof

\subsection{Fibrations over K3 surfaces}
\label{_fibra_ove_K3_Subsection_}

The purpose of this subsection is to prove 
\ref{_subvarie_of_G^i_Proposition_}.
Consider the fibration $G^m(M)\stackrel p \arrow M$ over the 
K3 surface $M$. Recall that $G^m(M)$ was defined as
a fibration with fibers classifying 
the codimension-$m$ ideals in $J^m(M)$,
where $J^m(M)$ is the bundle of rings 
$J^m(M)\restrict x = \calo_M/{\goth m_x}^m$. 
There is a decreasing filtration
\begin{equation}\label{_filtra_on_J_Equation_} 
   J^i(M) \supset {\goth W}(M) \supset {\goth W}^2(M) \supset ..., 
\end{equation}
with 
\[ {\goth W}^i(M)\restrict x = {\goth m_x}^i \cdot
   \calo_M/{\goth m_x}^m
\]
Consider the bundle
$V = {\goth W}^{l-1}(M)/{\goth W}^{l}(M)$.

\hfill

\lemma\label{_V_stable_Lemma_}
Let $M$ be a complex K3 surface 
which is generic with respect to some hyperka\"hler
structure. Let $V$ the holomorphic vector bundle defined above,
$V= {\goth W}^{l-1}(M)/{\goth W}^{l}(M)$. Then $V$
is isomorphic to a symmetric power of the cotangent bundle of $M$.
Moreover, $V$ is (Mumford-Takemoto)
stable for all K\"ahler structures on $M$, and has no
proper subbundles.

{\bf Proof:} The first assertion is clear.
Let us prove stability of $V$.
{}From Yau's proof of Calabi conjecture, it follows
that for all K\"ahler classes on $M$, $M$ is equipped with 
the hyperk\"ahler metric in the same K\"ahler class.
The Levi-Civita connection on the
cotangent bundle $\Lambda^1(M)$ of a hyperk\"ahler manifold
$M$ is hyperholomorphic
(\ref{_hyperho_conne_Definition_}), and hence Yang-Mills 
(\ref{Yang-Mills_Definition_}). 
Therefore, $\Lambda^1(M)$ is stable, and $V$ polystable (a tensor power of
a Yang-Mills bundle is again Yang-Mills).
The holonomy group of $\Lambda^1(M)$ is obviously isomorphic to 
$SU(2)$. Therefore, the holonomy group of
$V= S^l(\Lambda^1(M))$ is also $SU(2)$. The symmetric
power of the tautological representation of 
$SU(2)$ is obviously irreducible. Therefore, the
holonomy representation of $V$ is irreducible,
and $V$ cannot be represented as a direct sum of 
vector bundles. To show that $V$ has no proper subbundles,
we notice that $H^{1,1}(M) \cap H^2(M, \Z)=0$ because
$M$ is generic with respect to $\c H$. Therefore,
all coherent sheaves on $M$ have first Chern class zero.
We obtain that a proper subbundle of $V$ is destabilizing,
which contradicts stability of $V$.
\endproof

\hfill

Let $J^m_{gr}(M)$ be the graded sheaf of rings
associated with the filtration \eqref{_filtra_on_J_Equation_}.
Conside the fibration $G_{gr}^m(M)$ with the points classifying
codimension-$m$ ideals in the fibers of $J^m_{gr}(M)$.
There is a natural map $G^m(M)\stackrel \phi\arrow G_{gr}^m(M)$
associating to an ideal its associated graded quotient. 
Composing $\phi$ with the map $s:\; M \arrow G^{\frac{l(l-1)}{2}}(M)$,
we obtain the section $s_{gr}:\; M \arrow G_{gr}^{\frac{l(l-1)}{2}}(M)$
of the natural projection $p_{gr}:\; G_{gr}^{\frac{l(l-1)}{2}}(M)\arrow M$.
The following \ref{_subvarie_of_G^i_gr_Proposition_} obviously 
implies \ref{_subvarie_of_G^i_Proposition_}.

\hfill

\proposition\label{_subvarie_of_G^i_gr_Proposition_}
Let $M$ be a complex K3 surface.
Assume that $M$ admits a hyperk\"ahler structure $\c H$
such that $M$ is generic with respect to $\c H$.
Let $Y\subset G^m_{gr}(M)$ be
a closed irreducible subvariety of the total space of the
fibration  $G^m_{gr}(M) \stackrel p \arrow M$. 
Assume that $Y$ is generically finite over $M$.
Then $m = \frac{l\cdot (l+1)}{2}$ for some positive 
integer $l$, and $Y$ is the image of the natural map
$s_{gr}:\; M \arrow G^m_{gr}(M)$ constructed above.

\hfill

{\bf Proof:}
By \ref{_V_stable_Lemma_}, the bundle
$J^m_{gr}(M)$ is polystable. As usually, applying the Uhle\-n\-beck-\--Yau
theorem (\ref{_UY_Theorem_}), we endow the fibration 
$G^m_{gr}(M) \stackrel p \arrow M$ with a natural connection $\nabla$.
{}From \ref{_subvarie_of_PV_Proposition_} it is easy to deduce that
the image of $\pi:\; X \arrow G^m_{gr}(M)$ is preserved by 
the connection $\nabla$. Since  
$Y$ is generically finite over $M$,
the natural projection $Y\stackrel p \arrow M$ 
is a finite covering. Since $M$ is simply connected,
this map is an isomorphism. 

For $x\in M$, let $t_x \in X$ be the ideal of $J_{gr}^m(M)$
such that $p(t_x) = x$. Denote by $l$ the maximal
number such that $t_x \not\supset {\goth W}^{l-1}(M)$ for
some $x$. Consider the space 
\begin{multline*} A_x:= {\goth W}^{l-1}(M)\cap 
t_x \bigg/{\goth W}^{l}(M)\  \ \ \text{\LARGE $\subset$} \ \ \  {\goth W}^{l-1}(M)
   \bigg/ {\goth W}^{l}(M).\\
\end{multline*}
Let $w= \dim A_x$. Since $t_x$ is preserved by the connection
$\nabla$, the number $w$ does not depend on $x$.
This gives a $w$-dimensional subbundle $A$ in
$V= {\goth W}^{l}(M)\subset {\goth W}^{l-1}(M)$.
By \ref{_V_stable_Lemma_}, $A$ is either $V$
or empty. Since $t_x \not\supset {\goth W}^{l-1}(M)$,
$A=0$. Since $t_x$ is an ideal, this implies that
$t_x \subset {\goth W}^{l}(M)$. By definition
of $l$, it is the maximal number for which
$t_x \not\supset {\goth W}^{l-1}(M)$, and thus,
$t_x \supset {\goth W}^{l}(M)$. Therefore,
$t_x=l$. This proves
\ref{_subvarie_of_G^i_gr_Proposition_}.
We finished the proof of \ref{_subvarie_of_G^i_Proposition_}
and \ref{_gene_fini_universa_Theorem_}. \endproof


\section{Special subvarieties of the Hilbert scheme}
\label{_specia_subva_Section_}


\subsection{Special subvarieties}
\label{_special_subva_Subsection_}

\definition\label{_specia_subva_Definition_}
(See also \ref{_Unive_subva_Definition_}).
Let $M$ be a complex surface, $\c A$ a finite set,
$\phi:\; \c A \arrow M$ an arbitrary map.
For $i\in \c A$, consider the local ring $\calo_{\phi(i)}$
of germs of holomorphic functions in $\phi(i)$.
For $U\subset M$, consider the set $A_U$ of all
automorphisms (global or infinitesimal) of $U$ which
fix the image $\im\phi\subset M$ and act trivially on
$\calo_{\phi(i)}$. For $\gamma\in A_U$, we denote
by $\gamma^{[n]}$ the corresponding automorphism
of the Hilbert scheme $U^{[n]}$. A closed subvariety
$X\subset M^{[n]}$ is called {\bf special}
if for all $U\subset M$, all $\gamma\in A_U$,
$X\cap U^{[n]} $ is fixed by  $\gamma^{[n]}$.

\hfill

We are going to characterize special subvarieties
more explicitly, in the spirit of \ref{_inva_subva_from_Young_Theorem_}.

\hfill

Let $M$ be a complex surface, $\alpha$ a Young diagram.
\[ 
   \alpha = (n_1 \geq n_2 ...\geq n_k), \ \ \sum n_i = n,
\]
$\c A\subset \{1, ... , k\}$, and $\phi:\; \c A \arrow M$
an arbitrary map.
Consider a subvariety 
$\Delta_{(\alpha)}(\c A, \phi)$ of $M^{(n)}$
defined as in \eqref{_Delta(A,phi)_alpha_definition_Equation_}. 
A generic point $a \in \Delta_{(\alpha)}(\c A, \phi)$
is the one satisfying
\begin{equation}\label{_gene_of_Delta_alpha(A,phi)_Equation_}
\begin{split}
a = &\sigma(x_1, ..., x_n),\ \ \text{where} \ \  
         x_1 = x_2 = ... = x_{n_1}  \\&
               x_{n_1+1} = x_{n_1+2} = ... = x_{n_1+n_2} \  
               ..., \\ &
               x_{\sum_{i=1}^{k-1} n_i+1} = 
               x_{\sum_{i=1}^{k-1} n_i+2} = ... = x_n \\[4mm]
& 
\begin{minipage}[t]{0.7\linewidth}
and
\begin{description}
\item[(i)] $x_i = \phi(i)$ for all $i\in \c A$
\item[(ii)]
the points $x_1$, $x_{n_1+1}$, ..., 
$x_{\left(\sum_{i=1}^{j-1} n_i\right)+1}, ..., $ are pairwise unequal
for all $j\notin \c A$
\item[(iii)] the points $x_1$, $x_{n_1+1}$, ..., 
$x_{\left(\sum_{i=1}^{j-1} n_i\right)+1}, ..., $ don't belong 
to the set $\phi(\c A)$, for all $j\notin \c A$
\end{description}
\end{minipage}
\end{split}
\end{equation}

\hfill

We split the Young diagram 
\[ 
   \alpha = (n_1 \geq n_2 ...\geq n_k), \ \ \sum n_i = n,
\]
onto two diagrams,
$
   \alpha_{\c A} = (n_{a_1} \geq n_{a_2} ...\geq n_{a_l}), 
$
with $a_i$ running through $\c A$, and
$
   \check\alpha_{\c A} = (n_{b_1} \geq n_{b_2} ...\geq n_{b_{k-l}}), 
$
where $b_i$ runs through $\{1, ... , k\} \backslash \c A$.
Consider the Hilbert scheme $M^{[n_a]}$, where
$n_a = \sum n_{a_i}$. The map $\phi:\; {\c A} \arrow M$
gives a point $\Phi\in \Delta_{(\alpha_{\c A})}\subset M^{(n_a)}$
(see Subsection \ref{_speci_of_K3_Subsection_} for details).
Let $F_\phi:= \pi^{-1}(\Phi)$ be the fiber
of the standard projection $\pi:\; M^{[n_a]}\arrow M^{(n_a)}$.
For $y\in \Delta_{(\check\alpha_{\c A})}$ a generic point of
$\Delta_{(\check\alpha_{\c A})}$, let $F_{\check\alpha_{\c A}}(y)$
be the fiber of $\pi:\; M^{[n_b]}\arrow M^{(n_b)}$ over $y$.
Clearly, for $z\in \Delta_{(\alpha)}(\c A, \phi)$
a generic point, the fiber of $\pi:\; M^{[n]}\arrow M^{(n)}$
over $z$ is isomorphic to $F_\phi\times F_{\check\alpha_{\c A}}(y)$.
This isomorphism is not canonical, but is defined up to a twist by
the action of the group 
$G_y$ (see \ref{_inva_subva_F_alpha_identifi_Lemma_} for details).

\hfill

Fix a $G_y$-invariant
subvariety $E\subset F_\phi\times F_{\check\alpha_{\c A}}(y)$. 
For a generic point 
$z\in \Delta_{(\alpha)}(\c A, \phi)$, consider
a subvariety $E_z\subset \pi^{-1}(z)\subset M^{[n]}$
corresponding to $E$
under the isomorphism
\begin{equation} \label{_fiber_gene_special_subva_Equation_}
  \pi^{-1}(z) 
  \cong F_\phi\times F_{\check\alpha_{\c A}}(y).
\end{equation}
Let $\Delta_{[\alpha]}(\c A, \phi, E)$
be the closure of the union of $E_z$ for all 
$z\in\Delta_{(\alpha)}(\c A, \phi)$ 
satisfying \eqref{_gene_of_Delta_alpha(A,phi)_Equation_}.
Clearly, $\Delta_{[\alpha]}(\c A, \phi, E)$ 
is a closed subvariety in $M^{[n]}$.

\hfill

\theorem\label{_special_subva_explici_Theorem_}
Let $M$ be a complex surface, $\c A$ a finite set,
$\phi:\; \c A \arrow M$ an arbitrary map, and
$X\subset M^{[n]}$ a special subvariety, associated with $\phi$.
Then 
\begin{description}
\item[(i)]
there exist a Young diagram $\alpha$
\[ 
   \alpha = (n_1 \geq n_2 ...\geq n_k), \ \ \sum n_i = n,
\]
an injection $\c A \hookrightarrow 
\{1, ... , k\}$, and a $G_y$-invariant 
subvariety \[ E\subset F_\phi\times F_{\check\alpha_{\c A}}(y),\] 
such that  $X=\Delta_{[\alpha]}(\c A, \phi, E)$,  where
$F_\phi\times F_{\check\alpha_{\c A}}(y)$ 
and $\Delta_{[\alpha]}(\c A, \phi, E)$
are varieties constructed above. 
\item[(ii)]
Conversely,
$\Delta_{(\alpha)}(\c A, \phi, E)$ is a special
subvariety of $M^{[n]}$ for all $\c A$, $\pi$, $E$.
\end{description}

\hfill

{\bf Proof:} We use the notation of \ref{_specia_subva_Definition_}.
For sufficiently small $U$, the automorphisms
from $A_{U\backslash \im \phi}$ act $n$-transitively
on $U\backslash \im \phi$. This implies that
$\pi(X) = \Delta_{(\alpha)}(\c A, \phi)$,
for an appropriate Young diagram 
\[ 
   \alpha = (n_1 \geq n_2 ...\geq n_k), \ \ \sum n_i = n,
\]
and an embedding $\c A \hookrightarrow \{1, 2, ..., k\}$.
Let $x$ be a generic point of $\Delta_{(\alpha)}(\c A, \phi)$.
 Consider the isomorphism
$\pi^{-1}(x) \cong F_\phi\times F_{\check\alpha_{\c A}}(y)$
of \eqref{_fiber_gene_special_subva_Equation_}, and the
action of $G_y$ on $F_\phi\times F_{\check\alpha_{\c A}}(y)$. 
Clearly, $A_U$ acts on $\pi^{-1}(x)$ as $G_y$. Therefore,
the intersection $E_x:= X \cap \pi^{-1}(x)$ is $G_y$-invariant.
We intend to show that $X=\Delta_{[\alpha]}(\c A, \phi, E_x)$

For $x, y\in M^{(n)}$ generic points of $\Delta_{(\alpha)}(\c A, \phi)$,
there exists $U\supset \{x,y\}$ and an automorphism 
$\gamma:\; U \arrow U$ such that $\gamma^{(n)}(x) =y$,
for $\gamma^{(n)}:\; U^{(n)} \arrow U^{(n)}$ the induced by
$\gamma$ automorphism of $U^{(n)}$. Since $X$ is a special
subvariety, $\gamma^{[n]}$
maps $E_x$ to $E_y:= X \cap \pi^{-1}(y)$. By definition,
$\Delta_{[\alpha]}(\c A, \phi, E_x)$ is a closure of the
union of all $\gamma^{[n]}(E_x)$, for all $U\subset M$ and
$\gamma \in A_U$. On the other hand,
$X$ is a closure of the union of all $E_y$,
where $y$ runs through all generic
points of $\Delta_{(\alpha)}(\c A, \phi)$. Thus,
$X$ and $\Delta_{[\alpha]}(\c A, \phi, E_x)$ coinside.
This proves \ref{_special_subva_explici_Theorem_} (i).
\ref{_special_subva_explici_Theorem_} (ii) is clear.
\endproof

\subsection{Special subvarieties of the Hilbert scheme of K3}
\label{_speci_of_K3_Subsection_}

\theorem \label{_all_subva_are_special_Theorem_}
Let $M$ be a complex K3 surface admitting a hyperk\"ahler structure
$\c H$ such that $M$ is generic with respect to $\c H$,
$M^{[n]}$ its Hilbert scheme and $M^{(n)}$ its symmetric
power. Let $X\subset M^{[n]}$ be a closed irreducible subvariety
such that $X$ is generically finite over
$\pi(X) \subset M^{(n)}$. Assume that
$M$ has no holomorphic automorphisms. Then
$X$ is a special subvariety of $M^{[n]}$,
in the sense of \ref{_specia_subva_Definition_}. 

\hfill

{\bf Proof:} 
{}From \ref{_subva_in_M^(n)_Proposition_} it follows that
$\pi(X) = \Delta_{(\alpha)}(\c A, \phi)$ for appropriate
$\c A$, $\alpha$ and $\phi$.
As previously, we split the Young diagram 
$\alpha$ onto 
$
   \alpha_{\c A} = (n_{a_1} \geq n_{a_2} ...\geq n_{a_l}), 
$
with $a_i$ running through $\c A$, and
$
   \check\alpha_{\c A} = (n_{b_1} \geq n_{b_2} ...\geq n_{b_{k-l}}), 
$
where $b_i$ runs through $\{1, ... , k\} \backslash \c A$.
Let $n_a:= \sum n_{a_i}$, $n_b:= \sum n_{b_i}$.
Consider the natural map 
\begin{equation}\label{_M_na_prod_M_nb_Equation_}
M^{(n_a)} \times M^{(n_b)} \stackrel s \arrow M^{(n)},
\end{equation}
defined in such a way as that to
map $\Delta_{(\alpha_{\c A})}\times 
\Delta_{(\check\alpha_{\c A})}$ to $\Delta_{(\alpha)}$.
This map is obviously finite.
Let $x=(x_1, ... , x_{n_a})\in \Delta_{(\alpha_{\c A})}$
$y=(y_1, ... , y_{n_a})\in \Delta_{(\check\alpha_{\c A})}$
be the points satisfying $x_i \neq y_j$ $\forall i, j$.
Then the fiber of $\pi:\; M^{[n]} \arrow M^{(n)}$
in $s(x,y)$ is naturally isomorphic to the product
$\pi^{-1}(x) \times \pi^{-1}(y)$, where the first
$\pi$ is the standard projection $\pi:\; M^{[n_a]} \arrow M^{(n_a)}$
and the second one is the standard projection 
$\pi:\; M^{[n_b]} \arrow M^{(n_b)}$. Denote thus
obtained map
\begin{equation}\label{_M_na_prod_M_nb_fibers_Equation_}
\pi^{-1}(x) \times \pi^{-1}(y) \tilde\arrow \pi^{-1}(s(x,y))
\end{equation}
by $\theta$.
Together,
the maps \eqref{_M_na_prod_M_nb_Equation_}, 
\eqref{_M_na_prod_M_nb_fibers_Equation_}
give a correspondence
\[ \c D \subset \bigg( \Delta_{[\alpha_{\c A}]}\times 
                       \Delta_{[\check\alpha_{\c A}]}\bigg)
          \times\Delta_{[\alpha]}
\]
which is generically one-to-one over the first component
and generically finite over the second one. Denote the corresponding
projections from $\c D$ by $\pi_1$, $\pi_2$. Consider $X$
(the subvariety of $M^{[n]}$ given as data of 
\ref{_all_subva_are_special_Theorem_}) as a subvariety
of $\Delta_{[\alpha]}$. Let $\c D_X:=\pi_1(\pi_2^{-1}(X))$.
and $\Phi\in \Delta_{(\alpha_{\c A})}$ be the point given by
$\phi$, 
\[ 
   \Phi= \left(\underbrace{\phi(a_1), ..., \phi(a_1)}_
                                {n_{a_1} \text{\ times}}
   , \ \underbrace{\phi(a_1), ..., \phi(a_2)}_
                        {n_{a_2} \text{\ times}}, ...
   \right).
\]
Let $p_1$, $p_2$ be the 
projections of $\Delta_{[\alpha_{\c A}]}\times 
                       \Delta_{[\check\alpha_{\c A}]}$
to its components.
 Since
$\pi(X) =\Delta_{(\alpha)}(\c A, \phi)$, and $X$ is generically
finite over $\pi(X)$, the subvariety
\[ D_X\subset \Delta_{[\alpha_{\c A}]}\times 
\Delta_{[\check\alpha_{\c A}]}
\] 
is generically finite over 
$\{\Phi\}\times\Delta_{(\check\alpha_{\c A})}$.
Therefore, $p_2(\c D_2)\subset \Delta_{[\check\alpha_{\c A}]}$ is 
generically finite over $\{\Phi\}\times\Delta_{(\check\alpha_{\c A})}$.
Applying \ref{_gene_fini_universa_Theorem_}, 
we obtain that $p_2(\c D_X)$ is a universal
subvariety of $\Delta_{[\check\alpha_{\c A}]}$
The varieties $\Delta_{[\alpha_{\c A}]}$,
                       $\Delta_{[\check\alpha_{\c A}]}$
are equipped with the local action of the automorphisms $A_U$
(see \ref{_specia_subva_Definition_}).
Since $p_2(\c D_X)$ is universal,
\[ \c D_X\subset\{\Phi\}\times 
                       \Delta_{[\check\alpha_{\c A}]}
          \subset \Delta_{[\alpha_{\c A}]}\times 
                       \Delta_{[\check\alpha_{\c A}]}
\]
is fixed by the $A_U$-action. Therefore, 
$p_2(\c D_2)\subset \Delta_{[\check\alpha_{\c A}]}$
is also fixed by $A_U$. By construction, $\pi_2(\c D_X)=X$,
and thus, $X$ is fixed by $A_U$, i. e., special.
\endproof

\subsection{Special subvarieties of relative dimension 0}

\definition
Let $M$ be a complex surface, $M^{[n]}$ its Hilbert
scheme and $X\subset M^{[n]}$ an irreducible special subvariety.
The {\bf relative dimension of $X$} is the dimension of
the generic fiber of the projection
$\pi\restrict X :\; X \arrow \pi(X)$, where $\pi:\; M^{[n]}\arrow M^{(n)}$
is the standard morphism. 

\hfill

Let $\Delta_{(\alpha)}(\c A, \phi)\subset M^{(n)}$ be the subvariety
defined as in Subsection \ref{_subva_of_symme_special_Subsection_}.
Split $\alpha$ onto $\alpha_{\c A}$ and $\check \alpha_{\c A}$,
as in Subsection \ref{_special_subva_Subsection_}:
$
   \alpha_{\c A} = (n_{a_1} \geq n_{a_2} ...\geq n_{a_r}), 
$
with $a_i$ running through $\c A$, and
$
   \check\alpha_{\c A} = (n_{b_1} \geq n_{b_2} ...\geq n_{b_{k-r}}), 
$
where $b_i$ runs through $\{1, ... , k\} \backslash \c A$.
Let $\Phi\in \Delta_{(\alpha_{\c A})}$ be the point defined in
Subsection \ref{_speci_of_K3_Subsection_}. Consider the variety
$\pi^{-1}(\Phi)\subset M^{[n_a]}$, where $n_a = \sum n_{a_i}$. 

\hfill

\noindent
\proposition\label{_speci_unive_dime_zero_Proposition_}
\begin{description}
\item[(i)]
There exists a special subvariety $X\subset M^{[n]}$ such that
$\pi(X)= \Delta_{(\alpha)}(\c A, \phi)$ if and only if
all the numbers $n_{b_i}$ are of form $\frac{l\cdot (l+1)}{2}$, for 
integer $l$'s.

\item[(ii)]
 Let ${\goth S}(\alpha, \c A, \phi)$ be the set of all such subvarieties.
Assume that 
all the numbers $n_{b_i}$ are of form $\frac{l\cdot (l+1)}{2}$, for 
integer $l$'s.
Then ${\goth S}(\alpha, \c A, \phi)$ is in bijective 
correspondence with the set of points of 
$\pi^{-1}(\Phi)\subset M^{[n_a]}$.
\end{description}

{\bf Proof:} 
Using notation of the proof of 
\ref{_all_subva_are_special_Theorem_}, we
consider the subvariety 
$p_2(\c D_X)\subset \Delta_{(\check\alpha_{\c A})}$.
We have shown that this subvariety is 
universal of relative dimension 0. 
Therefore, ${\goth S}(\alpha, \c A, \phi)$
is nonepmpty if and only if all $n_{b_i}$ are of the 
form $\frac{l\cdot (l+1)}{2}$, for integer $l$'s.
Assume that ${\goth S}(\alpha, \c A, \phi)$ is nonempty.
Consider the unique universal subvariety 
$S\subset \Delta_{[\check\alpha_{\c A}]}$
of relative dimension 0,
constructed in \ref{_unive_subva_rela_dime_0_Proposition_}.
As in \ref{_inva_subva_from_Young_Theorem_}, a universal
subvariety of $\Delta_{[\check\alpha_{\c A}]}$ corresponds
to a $G_y$-invariant subvariety of the general fiber
of the projection
$\Delta_{[\check\alpha_{\c A}]} \arrow \Delta_{(\check\alpha_{\c A})}$.
Since $S$ is of relative dimension 0, the corresponding
$G_y$-invariant subvariety is a point. Denote this
point by $s$. Choose a point $f\in \pi^{-1}(\Phi)$.
Using the notation of \ref{_special_subva_explici_Theorem_},
and an isomorphism \eqref{_M_na_prod_M_nb_fibers_Equation_},
we construct a special subvariety 
$X\subset M^{[n]}$, 
$X = \Delta_{[\alpha]}(\c A, \phi, \{s\}\times \{f\})$.
{}From \ref{_special_subva_explici_Theorem_} it follows
that all special subvarieties of relative dimension
0 are obtained this way. Since $s$ is defined canonically, the
only freedom of choice we have after 
$\alpha, \c A, \phi$ are fixed 
is the choice of $f\in \pi^{-1}(\Phi)$.
This finishes the proof of 
\ref{_speci_unive_dime_zero_Proposition_}. \endproof

\hfill

\proposition \label{_symple_subva_special_Proposition_}
Let $M$ be a complex K3 surface with no complex
automorphisms, $M^{[n]}$ its Hilbert scheme
and $\Omega$ be the canonical
holomorphic symplectic form on $M^{[n]}$.
Assume that $M$ admits a hyperk\"ahler
structure $\c H$ such that $M$ is Mumford-Tate generic
with respect to $\c H$.  Let $X$ be an irreducible
complex subvariety of $M^n$, such that the
restriction $\Omega\restrict X$ is non-degenerate somewhere
in $X$.\footnote{Clearly, for $X$ trianalytic, $X_{ns}$ the non-singular
part of $X$, $\Omega\restrict{X_{ns}}$ is nowhere degenerate.}
 Then $X\subset M^{[n]}$ is a special 
subvariety of relative dimension 0.

{\bf Proof:} By 
\ref{_all_subva_are_special_Theorem_}, to prove that
$X$ is special it suffices
to show that $X$ is generically finite over
$\pi(X)$. This follows \ref{_triana_finite_in_gene_Claim_}.
\endproof

\hfill

\corollary\label{_one-to-one_triana_Corollary_}
Let $M$ be a complex K3 surface which is Mumford-Tate 
generic with respect to some hyperka\"hler structure,
$M^{[n]}$ its Hilbert scheme and $M^{(n)}$ its symmetric
power. Assume that $M$ has no holomorphic automorphisms.
Consider an arbitrary hyperka\"hler structure on $M^{[n]}$
whcih is compatible with the complex structure.
Let $X\subset M^{[n]}$ be a trianalytic subvariety
of $M^{[n]}$. Then $X$ is generically one-to-one
over $\pi(X)\subset M^{(n)}$ 

{\bf Proof:}
By \ref{_symple_subva_special_Proposition_},
$X$ is a special subvariety of relative dimension 0.
Now \ref{_one-to-one_triana_Corollary_}
follows from an explicit description of
 special subvarieties of relative dimension 0, given in
the proof of \ref{_speci_unive_dime_zero_Proposition_}.
\endproof


\section[Trianalytic and universal 
subvarieties of the Hilbert scheme of 
a general K3 surface]
{Trianalytic and universal 
subvarieties of the \\Hilbert scheme of 
a general K3 surface}
\label{_triana_unive_subva_Section_}


The aim of this section is to show that
 all trianalytic subvarieties of the Hilbert scheme
of a generic K3 surface are universal 
(\ref{_triana_subva_universal_Theorem_}).

\subsection{Deformations of trianalytic subvarieties}

We need the following general results on the structure of deformations
of trianalytic subvarieties, proven in \cite{_Verbitsky:Deforma_}. 
\footnote{The Desingularization Theorem
(\ref{_desingu_Theorem_}) significantly simplifies some
of the proofs of \cite{_Verbitsky:Deforma_}. 
This simplification is straightforward.} 

\hfill

Let $M$ be a compact hyperk\"ahler manifold,
$I$ an induced complex structure and $X\subset M$
a trianalytic subvariety. Consider $(X, I)$
as a closed complex subvariety of $(M,I)$. Let
$\Def_I(X)$ be a Douady space of $(X,I) \subset (M,I)$, that is,
a space of deformations of $(X,I)$ inside of $(M,I)$. 
By \ref{_G_M_invariant_implies_trianalytic_Theorem_}, 
for $X'\subset (M,I)$ a complex deformation
of $X$, the subvariety $X'\subset M$ is trianalytic. 
In particular, $X'$ is equipped with a natural singular hyperk\"ahler
structure (\cite{_Verbitsky:Deforma_}), 
i. e., with a metric and a compatible
quaternionic structure.

\hfill

\noindent
\theorem\label{_deforma_triana_Theorem_}
\begin{description}
\item[(i)] The $\Def_I(X)$ is a singular hyperk\"ahler variety,
which is independent from the choice of an induced complex 
structure $I$.
\item[(ii)] Consider the universal family $\pi:\; \c X \arrow \Def_I(X)$ 
of subvarieties of $(M,I)$, parametrized by $\Def_I(X)$.
Then the fibers of $\pi$ are isomorphic as hyperk\"ahler varieties.
\item[(iii)] Applying the desingularization functor
to $\pi:\; \c X \arrow \Def_I(X)$, we obtain a projection
$\pi:\; \tilde X\times Y \arrow Y$, where $Y$ is a desingularization
of $\Def_I(X)$ and $\tilde X$ is a desingularization of $X$.
\item[(iv)] The variety $\tilde X\times Y$ is equipped with a natural
hyperk\"ahler immersion to $M$.
\end{description}

{\bf Proof:} \ref{_deforma_triana_Theorem_} (i) 
and (ii) is proven in \cite{_Verbitsky:Deforma_}, and
\ref{_deforma_triana_Theorem_} (iii) is a trivial
consequence of \ref{_deforma_triana_Theorem_} (ii) and
the functorial property of the hyperk\"ahler 
desingularization. To prove 
\ref{_deforma_triana_Theorem_} (iv), we notice that
$\c X$ is equipped with a natural
morphism $f:\; \c X \arrow M$, which is
compatible with the hyperk\"ahler structure. 
Let $n:\;\tilde X\times Y\arrow \c X$ be the desingularization 
map. Clearly, the composition 
$\tilde X\times Y\stackrel n \arrow \c X \stackrel f\arrow M$
is compatible with the hyperk\"ahler structure. A 
morphism compatible with a hyperk\"ahler structure
is necessarily an isometry, and an isometry is always
an immersion.
\endproof

\subsection{Deformations of trianalytic
special subvarieties}

Let $M$ be a K3 surface, without automorphisms,
which is Mumford-Tate generic with respect to 
some hyperk\"ahler structure. 
In this Subsection, we study the deformations
of the special subvarieties of $M^{[n]}$.

\hfill

By \ref{_unive_subva_rela_dime_0_Proposition_},
universal subvarieties of $M^{[n]}$ are rigid. For the special
subvarieties, a description of its deformations
is obtained as an easy consequence of
\ref{_speci_unive_dime_zero_Proposition_}.

\hfill

\claim
Let $X\subset M^{[n]}$ be a special subvariety of relative
dimension 0, 
\[ X = \Delta_{[\alpha]}(\c A, \phi, \{s\}\times \{\psi\})\]
associated with $\Delta_{(\alpha)}(\c A, \phi)$
and $\psi\in F_{\phi}$ as in 
\ref{_speci_unive_dime_zero_Proposition_}.
Then the deformations of $X$ are locally parametrized
by varying $\phi:\; \c A \arrow M$ and $\psi\in F_{\phi}$.

\endproof

\hfill

Let $a_1$, ... , $a_r$ enumerate $\c A \subset \{ 1, ... , k\}$. Unless all
$n_{a_i}=1$, the dimension of $F_\phi$ is non-zero. Thus, the
union $\c X$ of all deformation of 
$X= \Delta_{[\alpha]}(\c A, \phi, \{s\}\times \{\psi\})$
is not generically finite over $\pi(\c X)=\Delta_{[\alpha]}(\c A, \phi)$.
Together with \ref{_deforma_triana_Theorem_} and
\ref{_symple_subva_special_Proposition_}, this suggests
the following proposition.

\hfill

\proposition \label{_n_i_=1_for_i_in_A_Proposition_}
Let $M$ be a complex K3 surface with no automorphisms which 
is Mumford-Tate generic with respect to some hyperk\"ahler structure. 
Consider the Hilbert scheme $M^{[n]}$ as a complex manifold. 
Let $\c H$ be an arbitrary hyperk\"ahler structure on $M^{[n]}$
agreeing with this complex structure. Consider an irreducible
trianalytic subvariety $X\subset M^{[n]}$. By 
\ref{_symple_subva_special_Proposition_}, 
$X$ is a special subvariety of $M^{[n]}$,
$X= \Delta_{[\alpha]}(\c A, \phi, \{\psi\}\times \{s\})$.
Then $n_{i}=1$ for all $i\in \c A$.

\hfill

{\bf Proof:} Consider $X$ as a complex subvariety in the
complex variety $M^{[n]}$. The corresponding Douady space
is described by \ref{_deforma_triana_Theorem_}. 
Consider the diagonal $\Delta_{(\alpha)}\subset M^{(n)}$.
For a general point $a\in \Delta_{(\alpha)}\subset M^{(n)}$,
the fiber $\pi^{-1}(a)$ is naturally decomposed as
in \eqref{_M_na_prod_M_nb_fibers_Equation_}:
\[
    \pi^{-1}(x) \times \pi^{-1}(y) \tilde\arrow \pi^{-1}(r(x,y)),
\]
for $a = r(x,y)$, where 
$r:\; M^{(n_a)}\times M^{(n_b)}\arrow M^{(n)}$ is a morphism 
of \eqref{_M_na_prod_M_nb_Equation_}.
Consider the subvariety $\pi^{-1}(y)\times \{s\}\subset \pi^{-1}(a)$.
This subvariety is clearly $G_a$-invariant, and applying
\ref{_inva_subva_from_Young_Theorem_}, we obtain
a universal subvariety of $M^{[n]}$. Denote this universal
subvariety by $\Delta_{[\alpha]}(\c A)$.
Let $\c X$ be the union of all complex 
deformations of $X\subset M^{[n]}$. From 
\ref{_deforma_triana_Theorem_} 
it is clear that $\c X$ is trianalytic;
from \ref{_speci_unive_dime_zero_Proposition_}
it follows that $\c X = \Delta_{[\alpha]}(\c A)$.
By \ref{_deforma_triana_Theorem_}, 
$\c X$ is trianalytic in $M^{[n]}$.
{}From \ref{_symple_subva_special_Proposition_}
it follows that all trianalytic subvarieties
$\c X\subset M^{[n]}$ are generically finite
over $\pi(\c X)\subset M^{(n)}$.
Consider the Young diagram
$\alpha_{\c A}= (n_{a_1}\geq n_{a_2} \geq ...)$. The generic
fiber of thus obtained generically finite map 
\begin{equation} \label{_Delta_alpha(A)_to_Delta_Equation_}
   \pi:\; \Delta_{[\alpha)}(\c A)\arrow \Delta_{(\alpha)}
\end{equation}
is isomorphic to $\pi^{-1}(y)$, where $y$ is a generic
point of $\Delta_{(\alpha_{\c A})}$, and
$\pi$ a projection 
$\pi:\; \Delta_{[(\alpha_{\c A})]} \arrow \Delta_{(\alpha_{\c A})}$. 
The dimension of this fiber is equal to $\sum (n_{a_i}-1)$.
By construction, $a_i$ enumerates $\c A\subset \{1, ... , k \}$.
Thus, $n_{i}=0$ for all $i\in \c A$.
This proves \ref{_n_i_=1_for_i_in_A_Proposition_}.
\endproof

\subsection{Special subvarieties and holomorphic symplectic form}

The aim of the Subsection is the following statement.

\hfill

\proposition \label{_holo_symple_degene_on_speci_Proposition_}
Let $M$ be a complex surface equipped with a holomorphically
symplectic form, $M^{[n]}$ its $n$-th Hilbert scheme, and
\[ X=\Delta_{[\alpha]}(\c A, \phi, \{\psi\}\times \{s\}) \]
the special subvariety of relative dimension 0.
Consider the holomorphic symplectic form $\Omega$ on $M^{[n]}$.
Assume that for all $i\in \c A$, $n_i=1$. 
Assume, furthermore, that the normalizarion $\tilde X$ of $X$ is smooth,
and the pullback of $\Omega$ to $\tilde X$ is a nowhere degenerate
holomorphic symplectic form on $\tilde X$. 
Then $\c A$ is empty. 
\footnote{To say that $\c A$ is empty
is the same as to say that $X$ is a universal subvariety of $M$.}

\hfill

{\bf Proof:} Assume that $\c A$ is nonempty. By 
\ref{_G_M_invariant_implies_trianalytic_Theorem_},
all complex deformations of $X$ are trianalytic.
Clearly, $\Delta_{[\alpha]}(\c A, \phi', \{\psi\}\times \{s\})$
is a deformation of $X$ for every $\phi':\; \c A \arrow M$.
Therefore, we may assume that $\phi:\; \c A \arrow M$ is an 
embedding.

Let $x\in  \Delta_{(\alpha)}$ be an arbitrary point.
We represent $x$ as in \eqref{_Del-a_alpha_definition_Equation_}.
The {\bf $i$-th component} of $x$ is 
$x_{\left(\sum_{i=1}^{i} n_1\right)+1}$, 
in notation of \eqref{_Del-a_alpha_definition_Equation_}. 
The components are defined up to a permutation
$i \arrow j$, for $i$, $j$ satisfying $n_i=n_j$.
Fix $p, q\in \{1, ... , k\}$, $p\in \c A$, $q\notin \c A$. 
Let $\Pi_{pq}:\; X \arrow \{\phi(q)\} \times M$ be the map
associating to $x\in X$ the $n_p$-th and $n_q$-th
components of $\pi(x) \in \Delta_{(\alpha)} \subset M^{(n)}$.
Clearly, the map $\Pi_{pq}$ is correctly defined.
Let $\check \Pi_{pq}:\; X \arrow M^{(n-n_p-n_q)}$
be the map associating to $x\in X$ the rest of components 
$x_{\left(\sum_{i=1}^{i} n_1\right)+1}$, $i \neq p, q$
of $x$. Let $2^M$ be the set of subsets of $M$.
Consider the map $c:\; M^{(i)}\arrow 2^M$ associating
to $x\in M^{(i)}$ the corresponding subset of $M$.
For $t= (\phi(q), t_0)\in \{\phi(q)\} \times M$, denote by $c(t)$
the subset $\{ \phi(q), t_0\} \subset M$. Let $X_0\subset X$ be the 
set of all $x\in X$ such that $c(\Pi_{pq}(x))$ does not
intersect eith $c(\check\Pi_{pq}(x))$. 
A subvariety $C\subset X$ is called non-degenerately
symplectic if the holomorphic symplectic form on $X$ is 
nowhere degenerate on $C$. For any
$t\in\check \Pi_{pq}(X)\subset  M^{(n-n_p-n_q)}$,
the intersection $X_t :={\check \Pi_{pq}}^{-1}(t) \cap X_0$ is
smooth. The
holomorphically symplectic form in a tangent space
to a zero-dimensional
sheaf $S\in M^{[n]}$, $Sup(S) = A \coprod B$
can be computed separately for the part with support in $A$
and the part with support in $B$. 
Thus, $X_t$ must be non-degenerately symplectic.
On the other hand, $X_t={\check \Pi_{pq}}^{-1}(t) \cap X_0$ is easy to
describe explicitly. Let $M_0= M\backslash c(t)$,
where $c(t)$ is again $t$ considered as a subset on $M$.
Then $X_t$ is canonically isomorphic to a blow-up of $M_0$
in $\{\phi(q)\}$. This blow-up is obviously not non-degerenerately
symplectic. We obtained a contradiction. This concludes
the proof of \ref{_holo_symple_degene_on_speci_Proposition_}.
\endproof

\subsection{Applications for trianalytic subvarieties}

\ref{_holo_symple_degene_on_speci_Proposition_}
implies the following theorem, which is the main
result of this section.

\hfill

\theorem\label{_triana_subva_universal_Theorem_}
Let $M$ be a complex K3 surface with no automorphisms which 
is Mum\-ford-\-Tate generic with respect to some hyperk\"ahler structure. 
Consider the Hilbert scheme $M^{[n]}$ as a complex manifold. 
Let $\c H$ be an arbitrary hyperk\"ahler structure on $M^{[n]}$
agreeing with this complex structure. Consider a
trianalytic subvariety $X\subset M^{[n]}$. Then 
$X$ is a universal subvariety of $M^{[n]}$ of relative dimension
0.

{\bf Proof:} By \ref{_n_i_=1_for_i_in_A_Proposition_},
$X$ is a special subvariety of $M^{[n]}$, 
$X=\Delta_{[\alpha]}(\c A, \phi, \{\psi\}\times \{s\})$.
The $X$ is non-degenerately symplectic because it si trianalytic. 
Applying \ref{_holo_symple_degene_on_speci_Proposition_},
we obtain that $\c A$ is empty, and $X$ is universal
in $M^{[n]}$. \endproof

\hfill

\corollary 
In assumptions of \ref{_triana_subva_universal_Theorem_},
$\codim_\C X \geq 4$, unless $X= M^{[n]}$.

{\bf Proof:} \ref{_special_subva_explici_Theorem_}
classifies universal subvarieties of relative dimension
0. All such subvarieties correspond to
diagonals $\Delta_{(\alpha)}\subset M^{(n)}$,
with 
\[ 
   \alpha = (n_1 \geq n_2 ...\geq n_k), \ \ \sum n_i = n,
\]
with all $n_i$ of form $\frac{l(l+1)}{2}$, with integer $l$'s.
Thus, for $X\neq M^{[n]}$ we have $n_1\geq 3$.
On the other hand, $\codim_\C \Delta_{(\alpha)} = 2 \sum (n_i-1)$,
so $\codim_\C \Delta_{(\alpha)}\geq 4$. Finally, since
$X$ is of relative dimension $0$, 
$\dim X = \dim \Delta_{(\alpha)}$,
so $\codim_\C X =\codim_\C \Delta_{(\alpha)}\geq 4$.
\endproof


\section{Universal subvarieties of the Hilbert scheme and algebraic
properties of its cohomology}
\label{_last_Section_}


\subsection{Birational types of algebraic subvarieties of relative
dimension 0}

Let $M^{[n]}$ be a Hilbert scheme, $\alpha$ a Young diagram 
satisfying assumptions of \ref{_unive_subva_rela_dime_0_Proposition_}
and $\c X_\alpha\subset M^{[n]}$ the corresponding universal subvariety of
relative dimension 0. Consider the natural map 
$\pi:\; M^{[n]}\arrow M^{(n)}$ mapping the Hilbert scheme
to the symmetric product of $n$ copies of $M$. Clearly,
$\pi(\c X_\alpha)= \Delta_{(\alpha)}$, where 
$\Delta_{(\alpha)}$ is the stratum of $M^{(n)}$ corresponding
to the Young diagram $\alpha$ as in 
Subsection \ref{_subva_of_symme_special_Subsection_}.
Moreover, from the definition of $\c X_\alpha$
it is evident that $\pi:\; \c X_\alpha\arrow \Delta_{(\alpha)}$
is a birational isomorphism.

\hfill

Let 
\[
\alpha = \left(\vphantom{\sum}
  n_1 = , ... , n_{i_1}> n_{i_1+1}=, ... , =n_{i_1+i_2} > ,...,
  > n_{1+\sum_{j=1}^{l-1} i_j}= , ..., n_{\sum_{j=1}^{l} i_j}\right)
\]
be a diagram satisfying assumptions of 
\ref{_unive_subva_rela_dime_0_Proposition_}. Then $\Delta_{(\alpha)}$ 
is birationally isomorphic to the product 
$\prod_{j=1}^l M^{(i_j)}$. Thus, $\c X_\alpha$ is birational to
a hyperka\"hler manifold $\prod_{j=1}^l M^{[i_j]}$.

\hfill

\theorem\label{_hype_bira_Theorem_}
(Mukai)
Let $f:\; X_1 \arrow X_2$ be a birational isomorphism of
compact complex manifolds of hyperka\"hler type. 
Then the second cohomology of $X_1$ is naturally
isomorphic to the second cohomology of $X_2$,
and this isomorphism is compatible with the Hodge structure
and the Bogomolov-Beauville\footnote{For the definition and properties of 
Bogomolov-Beauville form, see Subsection \ref{_B_B_Hilbe_Subsection_}.} 
form on $H^2(X_1)$, $H^2(X_2)$.

{\bf Proof:} Well known;  see, e. g.
\cite{_Mukai:Sugaku_}, \cite{_Huybrechts_}. \endproof

\hfill

We obtain that $H^2(\c X_\alpha)$ is isomorphic to 
$\oplus_{j=1}^l H^2(M^{[i_j]})$. Therefore,
$\dim H^{2,0}(\c X_\alpha)>1$ unless $l=1$.
On the other hand, by Bogomolov's decomposition
theorem (\ref{_simple_mani_crite_Theorem_}),
a hyperka\"hler manifold with $H^1(X, \R)=0$,
$\dim H^{2,0}(X)>1$ is canonically isomorphic to a product
of two hyperka\"hler manifolds of positive dimension. 
We obtain the following result.

\hfill

\proposition \label{_triana_bira_M^(l)_Proposition_}
Let $M$ be a complex K3 surface with no complex
automorphisms, admitting  a hyperk\"ahler
structure $\c H$ such that $M$ is Mumford-Tate generic
with respect to $\c H$.
Let $\c X_\alpha$ be a trianalytic subvariety of $M^{[n]}$,
which is by \ref{_triana_subva_universal_Theorem_}
universal and corresponds to a Young diagram
$\alpha = (n_1\geq n_2\geq ...\geq n_l)$.
Assume that $\c X_\alpha$ is not isomorphic to a product
of two hyperka\"hler manifolds of positive dimension. 
Then $n_1 = n_2 =, ..., =n_l$, and $\c X_\alpha$
is birationally equivalent to $M^{[l]}$.

\endproof

\subsection {The Bogomolov-Beauville form on the Hilbert scheme}
\label{_B_B_Hilbe_Subsection_}

Let $X$ be a simple\footnote{See \ref{_simple_mani_crite_Theorem_}
for the definition of ``simple''.}
hyperka\"hler manifold. It is well known that $H^2(X)$ is equipped
with a natural non-degenerate symmetric pairing 
\[ (\cdot, \cdot)_{\c B}:\; H^2(X) \times H^2(X) \arrow \C \]
which is compatible with the Hodge structure and with 
the $SU(2)$-action. This pairing is defined up to a constant
multiplier, and it is a topological invariant of $X$.
For a formal definition and basic properties of this form,
see \cite{_Beauville_} (Remarques, p. 775), and also
\cite{_Verbitsky:cohomo_}, \cite{_coho_announce_}.

\hfill

For a Hilbert scheme $M^{[n]}$ of points on a K3 surface,
the form $(\cdot, \cdot)_{\c B}$ can be computed explicitly
as follows.

\hfill

Consider the map 
$\pi:\; M^{[n]}\arrow M^{(n)}$ from the Hilbert scheme
to the symmetric power of $M$. Clearly, the space $H^2(M^{(n)})$
is naturally isomorphic to $H^2(M)$. 
Let $\Delta_{n}\subset M^{[n]}$ be the singular locus of
the map $\pi$, and $[\Delta_{n}]\in H^2(M^{[n]})$ its fundamental
class. The following proposition gives a full description
of the Bogomolov-Beauville form on $H^2(M^{[n]})$ in terms of
the Poincare form on $H^2(M)$.

\hfill

\proposition\label{_B-B_on_H^2(M^[n])_Proposition_}
Let $M$ be a K3 surface, and $M^{[n]}$ its Hilbert scheme
of points. Consider the pullback map 
$\pi^*:\; H^2(M) = H^2(M^{(n)}) \arrow H^2(M^{[n]})$.
Then
\begin{description}
\item[(i)] The map $\pi^*:\; H^2(M) \arrow H^2(M^{[n]})$
is an embedding. We have a direct sum decomposition
\begin{equation} \label{_H^2_M^[n]_decompo_Equation_}
H^2(M^{[n]}) = \pi^*(H^2(M)) \oplus \C \cdot [\Delta_{n}],
\end{equation}
where $[\Delta_{n}]\in H^2(M^{[n]})$ is the cohomology class defined
above. 
\item[(ii)] The decomposition \eqref{_H^2_M^[n]_decompo_Equation_}
is orthogonal with respect to the Bogomolov-Beauville form 
$(\cdot, \cdot)_{\c B}$. The restriction of $(\cdot, \cdot)_{\c B}$ to
\[ H^2(M) = \pi^*(H^2(M))\subset H^2(M^{[n]})\]  is equal to
the Poincare form times constant.  

Since the form $(\cdot, \cdot)_{\c B}$ is defined up to a 
constant multiplier, we may assume that, after a rescaling,  
$(\cdot, \cdot)_{\c B}\restrict {H^2(M)}$ is equal
to the Poincare form.

\item[(iii)] After a rescaling required by (ii), we have
\[ ( [\Delta_{n}], [\Delta_{n}])_{\c B} = -2 (n-1). \]

\end{description}

{\bf Proof:} Well known; see, for instance, 
\cite{_Huybrechts:hyperkahle_} 2.2. \endproof

\hfill

Further on, we always normalize the Bo\-go\-mo\-lov-\-Beau\-ville form
on \linebreak $H^2(M^{[n]})$ as in \ref{_B-B_on_H^2(M^[n])_Proposition_}
(ii).

\subsection{Frobenius algebras associated with vector spaces}

Let $X$ be a compact hyperk\"ahler manifold.
The algebraic structure of $H^*(X)$ is studied using the
general theory of Lefschetz-Frobenius algebras, introduced
in \cite{_Lunts-Loo_}. For details of definitions and 
computations, the reader is referred to \cite{_Verbitsky:cohomo_},
\cite{_coho_announce_}.

\hfill

\definition
Let $A= \bigoplus\limits^{2d}_{i=0} A_i$
be a graded commutative associative algebra over a field of 
characteristic zero. 
Assume that $A_{2d}$ is 1-dimensional, and the natural 
linear form $\epsilon:\; A \arrow A_{2d}$ projecting $A$ to
a $A_{2d}$ gives a non-degenerate
scalar product $a, b \arrow \epsilon(ab)$. Then $A$ is called
{\bf a graded commutative Frobenius algebra}, or
Frobenius algebra for short.

\hfill

\proposition \label{_S^*_algebra_Fro_Proposition_}
Let $V$ be a vector space equipped with a non-\-de\-ge\-ne\-rate scalar 
product, and $n$ a positive integer number. Then there exist a unique up to
an isomorphism Frobenius algerba
\[ A(V, n) = A_0 \oplus A_2 \oplus .. \oplus A_{4n}
\]
such that
\begin{description}
\item[(i)] 
\begin{equation*}\begin{split}
A_{2i} = S^i(V), \ \ \ &\text{for $i\leq n$} \\
A_{2i} = S^{2n-i}(V), \ \ \ &\text{for $i\geq n$} 
\end{split}
\end{equation*}
and
\item[(ii)] For an operator $g\in SO(V)$,
consider the corresponding endomorphism of $S^*(V)$.
This way, $g$ might be considered as a linear operator
on $A$. Then $g$ is an algebra automorphism.
\end{description}
{\bf Proof:} \ref{_S^*_algebra_Fro_Proposition_} is elementary.
For a complete proof of existence and uniqueness
of $A(V, n)$, see \cite{_Verbitsky:cohomo_}. \endproof

\hfill

The importance of the algebra $A(V, n)$ is explained by the following
theorem.

\hfill

\theorem\label{_stru_of_H^*_Theorem_}
\cite{_Verbitsky:cohomo_}
Let $X$ be a compact connected simple hyperk\"ahler manifold.
Consider the space $V= H^2(X)$, equipped with the
natural scalar product of Bogomolov-Beauville
(Subsection \ref{_B_B_Hilbe_Subsection_}).
Let $A$ be a subalgebra of $H^*(X)$ generated by $H^2(M)$. 
Then $A$ is naturally isomorphic to $A(V,n)$.

\endproof

\subsection{A universal embedding from a K3 to its Hilbert scheme}

Consider a universal embedding 
$M \stackrel \phi \hookrightarrow M^{[n]}$,
$n = \frac{k (k-1)}{2}$, $n >1$, mapping a point
$x\in M$ to subscheme given by the ideal $({\goth m}_x)^k$,
where ${\goth m}_x$ is the maximal ideal of $x$.
Pick a hyperka\"hler structure $\c H$ on $M^{[n]}$.
The aim of this subsection is to prove the following
result.

\hfill

\proposition \label{_ima_of_M_not_triana_Proposition_}
The image of $\phi$ is not trianalytic in $M^{[n]}$.

\hfill

The proof of \ref{_ima_of_M_not_triana_Proposition_}
takes the rest of this section. Together with 
\ref{_unive_subva_rela_dime_0_Proposition_} and
\ref{_gene_fini_universa_Theorem_}, this result
immediately implies the following corollary.

\hfill

\corollary
Let $M$ be a complex K3 surface without automorphisms. 
Assume that $M$ admits a hyperk\"ahler
structure $\c H$ such that $M$ is Mumford-Tate
generic with respect to $\c H$ (\ref{_generic_manifolds_Definition_}).
Pick a hyperka\"hler structure $\c H'$ on its Hilbert scheme 
$M^{[n]}$ ($n>1$).
Let $X\subset M^{[n]}$ be a trianalytic subvariety
of $M^{[n]}$. Then $\dim_{\Bbb H} X>1$.

\endproof

\hfill

{\bf Proof of \ref{_ima_of_M_not_triana_Proposition_}:}
Consider the map \[ \phi^*:\; H^4(M^{[n]}) \arrow H^4(M) = \C. \]
To prove that $\im \phi$ is not trianalytic in $M^{[n]}$,
it suffices to show that $\phi^*$ is not $SU(2)$-invariant. 
Let $H^4_r(M^{[n]})$ be the subspace of 
$H^4(M^{[n]})$ generated by $H^2(M^{[n]})$,
and $f:\; H^4_r(M^{[n]})\arrow \C$ the restriction
of $\phi^*$ to $H^4_r(M^{[n]})\subset H^4(M^{[n]})$.
Since the subspace $H^4_r(M^{[n]})\subset H^4(M^{[n]})$
is $SU(2)$-invariant, the map $f$ should be $SU(2)$-invariant
if $\phi^*:\; H^4(M^{[n]}) \arrow H^4(M) = \C$ is $SU(2)$-invariant.

\hfill

By \ref{_stru_of_H^*_Theorem_}, the space $H^4_r(M^{[n]})$
is naturally isomorphic to \[ S^2H^2(M^{[n]}).\] Thus, $f$
can be considered as a map $f:\; S^2H^2(M^{[n]})\arrow \C$.
Consider the Bo\-go\-mo\-lov-\-Beau\-ville form as another such map
$B:\; S^2 H^2 (M^{[n]})\arrow \C$. From 
\ref{_B-B_on_H^2(M^[n])_Proposition_}, it is clear
that $f = B + 2 (n-1) d^2$, where $d:\; H^2(M^{[n]})\arrow C$
is the projection of $H^2(M^{[n]})$ to 
the component $\C = \C \cdot [\Delta_n]$ of the
decomposition \eqref{_H^2_M^[n]_decompo_Equation_}.
Since $B$ is $SU(2)$-invariant, the map $f$
is $SU(2)$-invariant if and only if the map
$d^2:\; S^2H^2(M^{[n]})\arrow \C$ is $SU(2)$-invariant.
Therefore, the following claim is sufficient to prove
\ref{_ima_of_M_not_triana_Proposition_}.

\hfill

\claim \label{_d^2_not_SU(2)_inv_Claim_}
In the above notations, the vector
$d^2\in S^2H^2(M^{[n]})^*$ is not $SU(2)$-invariant.

\hfill

{\bf Proof:}
Let $V$ be the $SU(2)$-subspace of $S^2H^2(M^{[n]})^*$ generated
by $d^2$. Acting on $d^2$ by various $g\in \goth{su}(2)$, we can obtain
any element of type $d \cdot g(d)$ (this follows from Leibnitz rule).
Therefore, $V = SU(2) \cdot d^2$ contains $d\otimes V_0$, where
$V_0 \subset H^2(M^{[n]})^*$ is the $SU(2)$-subspace
of $H^2(M^{[n]})^*$ generated by $d$. 
Acting on $d \cdot g(d)$ by various $h\in SU(2)$,
we obtain any element of type $h(d) \cdot h(gd)$.
This implies that $V = SU(2) \cdot d\otimes V_0$ contains
$S^2(V_0)$. We obtained that $V = S^2 V_0$.

Clearly, $d^2$ is 
$SU(2)$-invariant if and only if $V$ is 1-dimensional.
Thus, $d^2$ is $SU(2)$-invariant if and only if $V_0$ is 1-dimensional,
that is, if $d$ is $SU(2)$-invariant. On the other hand,
the map $d$ is an orthogonal projection to
$\C \cdot [\Delta_n]\subset H^2(M^{[n]})$. Thus,
$d$ is $SU(2)$-invariant if and only if 
$[\Delta_n]\in H^2(M^{[n]})$ is $SU(2)$-invariant. The class
$[\Delta_n]\in H^2(M^{[n]})$
is a fundamental class of a subvariety $\Delta_n \subset M^{[n]}$.
By \ref{_G_M_invariant_implies_trianalytic_Theorem_},
$[\Delta_n]$ is $SU(2)$-invariant if and only if 
$\Delta_n \subset M^{[n]}$ is trianalytic. The trianalytic
subvarieties are hyperka\"hler, outside of singularities.
Since $\Delta_n$ is a divisor, it has odd complex dimension
and cannot be hyperka\"hler. Thus, the class $[\Delta_n]$ is not 
$SU(2)$-invariant, and the map $d^2:\; S^2H^2(M^{[n]})\arrow \C$
is not $SU(2)$-invariant. This proves 
\ref{_d^2_not_SU(2)_inv_Claim_}. 
\ref{_ima_of_M_not_triana_Proposition_} is proven.
\endproof

\subsection{Universal subvarieties of the Hilbert scheme and 
Bo\-go\-mo\-lov-\-Beau\-ville form}

In this subsection, we show that the subvarieties 
$\c X_\alpha\subset M^{[n]}$, obtained as in
\ref{_triana_bira_M^(l)_Proposition_}, are not trianalytic.
We prove this using the explicit calculation of the 
Bogomolov-Beauville form on $M^{[i]}$ 
(\ref{_B-B_on_H^2(M^[n])_Proposition_})
and the following result.

\hfill

\claim \label{_pullback_B-B_SU(2)_inv_Proposition_}
Let $\phi:\; X \hookrightarrow Y$ be a morphism of compact
hyperka\"hler manifolds.
Consider the corresponding pullback map
\[ \phi^*:\; H^2(Y) \arrow H^2(X).\] Let 
\[ \Psi:\; S^2H^2(Y)  \arrow S^2H^2(X) \] be the symmetric
square of $\phi^*$, and $B_Y\in S^2H^2(Y)$ the vector
corresponding to the Bogomolov-Beauville pairing.
Then $\Psi(B_Y)$ is $SU(2)$-invariant, with respect to
the natural action of $SU(2)$ on $S^2 H^2(X)$.

\hfill

{\bf Proof:} It is well known that, for every 
morphism of hyperka\"hler varieties, the pullback map is compatible with the
$SU(2)$-action in the cohomology. To see this, one may notice
that the $SU(2)$-action is obtained from the Hodge-type grading
associated with induced complex structures, and the pullback
is compatible with the Hodge structure. Now, $B_Y$ is $SU(2)$-invariant,
and therefore, $\Psi(B_Y)$ is also $SU(2)$-invariant.
\endproof

\hfill

Let $M$ be a K3 surface, $M^{[n]}$ its Hilbert scheme and
$\c X_\alpha$ be a universal subvariety of $M^{[n]}$
of relative dimension 0, obtained from the Young diagram
$(n_1 = , ..., = n_l)$ as in \ref{_triana_bira_M^(l)_Proposition_}.
Assume that $\c X_\alpha$ is trianalytic with respect to some
hyperka\"hler structure on $M^{[n]}$. The manifold
 $\c X_\alpha$ is birational to $M^{[l]}$ 
(\ref{_triana_bira_M^(l)_Proposition_}).
By \ref{_hype_bira_Theorem_}, there is a natural
isomorphism $H^2(\c X_\alpha)\cong H^2(M^{[l]})$, 
and this isomorphism is compatible with the 
Bogomolov-Beauville form. Consider the map
$\phi^*:\; H^2(M^{[n]}) \arrow H^2(\c X_\alpha) = H^2(M^{[l]})$.
Recall that $H^2(M)$ is considered as a subspace
of $H^2(M^{[i]})$, for all $i$ (\ref{_B-B_on_H^2(M^[n])_Proposition_} (i))
Clearly, $\phi^*$ acts as identity on the subspaces
$H^2(M) \subset H^2(M^{[n]})$, $H^2(M) \subset H^2(M^{[l]})$.
Consider the pullback 
$\phi^*([\Delta_n])\in H^2(\c X_\alpha) = H^2(M^{[l]})$ 
of $[\Delta_n]\in H^2(M^{[n]})$. 

\hfill

\lemma\label{_pullba_of_Delta_Lemma_}
In the above notations, $\phi^*([\Delta_n]) = \frac{n}{l} [\Delta_l]$.

\hfill

{\bf Proof:} Let $t\in H^2(M) \subset H^2(M^{[l]})$
be the component of $\phi^*([\Delta_n])$ corresponding to the
decomposition \eqref{_H^2_M^[n]_decompo_Equation_}.
Since the decomposition \eqref{_H^2_M^[n]_decompo_Equation_}
is integer, $t$ is an integer cohomology class.
Let $\c M$ be a universal K3 surface, considered as
a fibration over the moduli space $D$ of marked K3 surfaces.
The construction of the Hilbert scheme can be applied
to the fibers of $\c M$. We obtain a universal Hilbert scheme
$\c M^{[n]}$, which is a fibration over $D$. 
Since $\c X_\alpha$ is a universal subvariety
of $M^{[n]}$, there is a corresponding fibration 
over $D$ as well. Consider the class $t\in H^2(M)$
as a function $t(I)$ of the complex structure $I$ on $M$.
Since the cohomology class $t(I)$ is integer, and $D$ is connected, 
$t(I)$ it is independent from $I\in D$. On the other hand,
$t(I)$ has type  $(1,1)$ with respect to $I$. There are no
non-zero cohomology classes $\eta\in H^2(M)$ which have
type $(1,1)$ with respect to all complex structures on $M$.
Thus, $t=0$. We obtain that $\phi^*([\Delta_n]) = r [\Delta_l]$,
where $r$ is some integer number. It remains to check that
$r= \frac{n}{l}$. Recall that $\frac{n}{l}$ is an integer number
which is equal to $\frac{k(k-1)}{2}$, for some $k\in \Z$,
$k>1$.

The points of the Hilbert scheme $M^{[i]}$ correspond to ideals
$I\in\calo_{M}$, $\dim \calo_{M}/I =i$. Consider a rational map
$\xi:\; M^{[l]} \arrow M^{[n]}$ mapping an ideal $I\subset \calo_{M}$
to $I^k$. Let $S\subset M^{[l]}$ be the union of all strata
$\Delta_{[\alpha]}$ of codimension more than 1.
It is easy to check that $\xi$ is well defined outside
of $S$: for all ideals $I\in M^{[l]}\backslash S$,
the ideal $I^k$ satisfies $\dim \calo_{M}/I^k =n$.
Consider the pullback map on the cohomology
associated with the morphism 
$\psi:\; M^{[l]}\backslash S \arrow M^{[n]}$.
Clearly, $H^2(M^{[l]}\backslash S)= H^2(M^{[l]})$.
The map $\phi^*:\; H^2(M^{[n]}) \arrow H^2(\c X_\alpha) = H^2(M^{[l]})$
is equal to 
\[ \xi^*:\; H^2(M^{[n]}) \arrow H^2(M^{[l]}\backslash S)=
   H^2(M^{[l]}) 
\]
Let $p\in \pi_2(M^{[l]}\backslash S)$ be an element of the second homotopy
group corresponding to $[\Delta_l]$ under Gurevich isomorphism.
It remains to show that $\xi(p)$ is equal to $\frac{n}{l}$ times
the element of the second homotopy group $\pi_2(M^{[n]})$
corresponding to $[\Delta_n]$.
The following observation is needed to understand the
geometry of $\Delta_i$.

\begin{equation}\label{_Arti_leng_2_Equation_}
\begin{minipage}[m]{0.8\linewidth}
Closed Artinian subschemes $\xi\subset M$ 
of length 2 with support in $x\in M$ are in one to one correspondence
with the vectors of projectivization of $T_x M$.
\end{minipage}
\end{equation}
Therefore, generic points of $\Delta_i$ correspond
to the triples $(\c X, x, \lambda)$, where
$\c X$ is a non-ordered set of $(i-1)$ distinct points of $M$,
$x\in M\backslash \c X$ a point of $M$
and $\lambda\in {\Bbb P} T_x M$ a line in $T_x M$.
Fix $\c X$, $x$ and an isomorphism ${\Bbb P} T_x M\cong \C P^1$.
We pick a map $p_x:\; S^2 \arrow M^{[l]}$ in such a way that
$p_x(\theta)= (\c X, x, \theta)$. Clearly, the corresponding
element of $\pi_2(M^{[l]})$ is mapped to $[\Delta_l]$ by
Gurevich's isomorphism.

To simplify notations, we
assume that $l=2$. It is easy to do the
case of general $l$ in the same spirit as we do $l=2$.

Let $U$ be a neighbourhood of $x$
in $M$. Taking $U$ sufficiently small, we may assume that
$U$ is equipped with coordinates. Let $R_\lambda$ be a 
parallel translation of $U$ along these coordinates
in the direction of $\lambda$, for
$\lambda \in \C^2$. The map $R_\lambda$ is defined in
a smaller neighbourhood $U'$ of $x$:
$R_\lambda:\; U_1 \arrow U$. The coordinates give
a natural identification $\C P^1 \cong {\Bbb P} T_y M$,
for all $y\in U$.
Let $\lambda\in T_x M\backslash 0$ be a vector
corresponding to $\theta\in {\Bbb P} T_x M$.
Clearly, then, 
\[ p_x(\theta) 
   = \lim\limits_{t\mapsto 0, t\in \R\backslash 0} \bigg\{x,
   R_{t\lambda}(x) \bigg\},
\]
where the pair $\{x,
   R_{t\lambda}(x) \}$ is considered as a point in $M^{[2]}$.

Let 
$y_i \in U_1^{[\frac{n}{2}]}\backslash \Delta_{\frac{n}{2}}$ 
be a sequence of points
converging to $\phi'(x)$, where $\phi':\; M \mapsto M^{[\frac{n}{2}]}$
maps a maximal ideal of $x\in M$ to its $k$-th power. 
Denote the support of $y_i$ by $S_{y_i}$.
Clearly, 
\[ \phi(p_x(\theta)) = \lim\limits_{t\mapsto 0, t\in \R\backslash 0} 
   \left(\vphantom{\prod}
   \lim \limits_{i\mapsto\infty}\{y_i, R_{t\lambda}(y_i) \}\right).
\] 
Taking limits in different order, we obtain
\[ \phi(p_x(\theta))= \lim \limits_{i\mapsto\infty}\left(\vphantom{\prod}
   \lim\limits_{t\mapsto 0, t\in \R\backslash 0} \{y_i,
   R_{t\lambda}(y_i) \} =  \lim \limits_{i\mapsto\infty} y_i(\theta)\right),
\]
where $y_i(\theta)\in \Delta_n$ is a point corresponding to 
a closed subscheme in $M^{[n]}$
with support $S_{y_i}$, of length 2 at every point 
of its support. For $y\in S_{y_i}$, consider the restriction
$y_i(\theta)\restrict y$
of $y_i(\theta)$ to $y$, which is a closed Artinian subscheme 
of length 2 with support in $y$. Clearly,
$y_i(\theta)\restrict y$
corresponds to $\theta\in {\Bbb P}T_y M$ as in
\eqref{_Arti_leng_2_Equation_}. Varying $\theta$, we obtain a
map $p_i:\; \C P^1 \arrow \Delta_n \backslash S$, 
$\theta \arrow y_i(\theta)$. By construction, this map is
homotopic to $\phi(p_x)$. On the other hand, it is clear
that $p_{y_i}$ is homotopic to $Card(S_{y_i})=\frac{n}{l}$ times a  
homotopy class represented by a map
$p_x:\; \C P^1 \arrow \Delta_n \backslash S$.
This proves \ref{_pullba_of_Delta_Lemma_}.
\endproof

\hfill

The Beauville-Bogomolov form identifies the second cohomology
with its dual. Thus, this form can be considered as a tensor
in the symmetric square of the second cohomology.
To show that the embedding $\c X_\alpha \hookrightarrow M^{[n]}$
is not trianalytic, we compute the pullback $\phi^* B_{M^{[n]}}$ of the
Beauville-Bogomolov tensor $B_{M^{[n]}}\in S^2 H^2(M^{[n]})$.
Consider the decomposition \eqref{_H^2_M^[n]_decompo_Equation_}
\begin{equation}\label{_H^2_M^[n]_decompo_again_Equation_}
H^2(M^{[i]}) = H^2(M) \oplus  \C \cdot \Delta_i. 
\end{equation}
Let $P\in S^2H^2(M)$ be the tensor corresponding to the Poincare
pairing. \ref{_B-B_on_H^2(M^[n])_Proposition_} computes the 
form $(\cdot,\cdot)_{\c B}\in S^2H^2(M)^*$ in terms of 
Poincare form and the decomposition 
\eqref{_H^2_M^[n]_decompo_again_Equation_}:
$(\cdot,\cdot)_{\c B}= P - 2(n-1)d^2$.  Therefore, the dual
tensor can be written as
$B_{M^{[i]}} = P - \frac{1}{2(n-1)}[\Delta_i]^2$.
We have shown that $\phi^*$ acts as an identity on the
first summand of \eqref{_H^2_M^[n]_decompo_again_Equation_},
and maps $[\Delta_n]$ to $\frac{n}{l}[\Delta_l]$.
Therefore,
\begin{equation} \label{_B_Ml_in_terms_of_phi_B_Mn_Equation_}
   \phi^* B_{M^{[n]}} = P - \frac{1}{2(n-1)}\frac{n}{l}[\Delta_l]^2 
   = B_{M^{[l]}} + 
   \left(\frac{1}{2(l-1)} -\frac{1}{2(n-1)}\frac{n}{l}\right) [\Delta_l]^2
\end{equation}
By \ref{_d^2_not_SU(2)_inv_Claim_}, $[\Delta_l]^2$ is not
$SU(2)$-invariant. Since $B_{M^{[l]}}$ {\it is} $SU(2)$-invariant,
$\phi^* B_{M^{[n]}}$ is $SU(2)$-invariant if and only if the coefficient
of $[\Delta_l]^2$ in \eqref{_B_Ml_in_terms_of_phi_B_Mn_Equation_}
vanishes:
\begin{equation} \label{_coeffi_pullback_formula_Equation_}
  (l-1)^{-1} -(n-1)^{-1}\frac{n}{l} =0 
\end{equation}
Clearly, this happens only if $(n- \frac{n}{l})= n-1$, i. e. when
$\frac{n}{l}=1$.
By definition, $\frac{n}{l} = \frac{k(k+1)}{2}\geq 3$. 
Therefore, $\phi^* B_{M^{[n]}}$ is not
$SU(2)$-invariant,
and $\c X_\alpha$ is not trianalytic in $M^{[n]}$.\footnote{
We do not need the whole strength of \ref{_pullba_of_Delta_Lemma_}
to show that the tensor $\phi^* B_{M^{[n]}}$ is not
$SU(2)$-invariant. Let $t$ be the coefficient of
\ref{_pullba_of_Delta_Lemma_}, $\phi^*([\Delta_n]) = t [\Delta_l]$.
To show that $\phi^* B_{M^{[n]}}$ is not
$SU(2)$-invariant, we need only to check that
$n-1 \neq t (l-1)$, where $n= l\frac{k(k+1)}{2}$.
Since $l\frac{k(k+1)}{2}-1$ is not divisible by 
$l-1$ for most $l$, $k$, the inequality $n-1 \neq t (l-1)$
holds automatically for most $l$, $n$.}

Comparing this with \ref{_pullback_B-B_SU(2)_inv_Proposition_}, 
and \ref{_triana_bira_M^(l)_Proposition_}, we obtain the following result.

\hfill

\theorem \label{_no_triana_subva_of_Hilb_Theorem_}
Let $M$ be a complex K3 surface without automorphisms. Assume that
$M$ is Mumford-Tate generic with respect to some hyperka\"hler structure. 
Consider the Hilbert scheme $M^{[n]}$ of points on $M$.
Pick a hyperk\"ahler structure on $M^{[n]}$ which is compatible with
the complex structure. Then $M^{[n]}$ has no proper
trianalytic subvarieties.

\endproof

\hfill

\remark 
It is easy to see that a generic K3 surface has no complex
automorphisms.

\hfill

\corollary \label{_no_comple_subva_of_gen_Hilb_Corollary_}
Let $M$ be a complex K3 surface.
Consider its Hilbert scheme $M^{[n]}$.
Let $\c M$ be the generic deformation of $M^{[n]}$.
Then $\c M$ has no complex subvarieties.

{\bf Proof:} Follows immediately from 
\ref{_no_triana_subva_of_Hilb_Theorem_}
and \ref{_hyperkae_embeddings_Corollary_}. 
\endproof

\hfill

{\bf Acknowledegments:} 
I am grateful to R. Bezrukavnikov and A. 
Beilinson for lecturing me on perverse sheaves and semismall resolutions, 
H. Nakajima for sending me the manuscript of \cite{_Nakajima_},
V. Lunts for fascinating discourse on the 
Mumford-Tate group, and to P. Deligne, M. Grinberg, 
D. Kaledin, D. Kazhdan and T. Pantev for valuable discussions.
My gratitude to  P. Deligne and D. Kaledin, who
found important errors in the earlier versions of this paper.

\end{document}